\documentclass[reqno]{amsart}
\usepackage{a4wide}
\usepackage{pstricks}
\usepackage{graphics,graphicx}
\usepackage[english]{babel}
\usepackage{amsfonts,amsthm,amssymb,amsmath,mathrsfs}

\usepackage[colorlinks=false, linktocpage=true]{hyperref}
\newtheorem{thm}{Theorem}[section]
\newtheorem{cor}[thm]{Corollary}
\newtheorem{lem}[thm]{Lemma}
\newtheorem{prop}[thm]{Proposition}
\theoremstyle{definition}
\newtheorem{defn}[thm]{Definition}

\theoremstyle{remark}
\newtheorem{rem}[thm]{Remark}
\numberwithin{equation}{section}

\usepackage[notcite,notref]{showkeys}

\begin{document}
\title[Pattern Formation for Reaction-Difussion Systems]{Two simple criterion to prove the existence of patterns in reaction-diffusion models of two components.}

\author{Francisco J. Vielma-Leal}\address{Departamento de Matemática, Facultad de Ciencias Naturales, Matemática y del Medio Ambiente, Universidad Tecnológica Metropolitana, Santiago, Chile.}
\email{fvielma@utem.cl}

\author{Miguel A.D.R. Palma}\address{IME, Universidade Federal Fluminense, Niter\'oi, Brazil.}
\email{mipalma@id.uff.br}

\author{Miguel Montenegro-Concha}\address{Departamento de Matemática, Facultad de Ciencias Naturales, Matemática y del Medio Ambiente, Universidad Tecnológica Metropolitana, Santiago, Chile.}
\email{miguel.montenegro@utem.cl}

\dedicatory{Dedicated to the memory of Marcos Lizana Peña}
\begin{abstract}
The aim of this work is to study the effect of diffusion on the stability of the equilibria in a general two-components reaction-diffusion system with Neumann boundary conditions in the space of continuous functions. As by product,  we establish sufficient conditions on the diffusive coefficients and other parameters for such a   reaction-diffusion model to exhibit patterns and we analyze their stability. We apply the results obtained in this paper to explore under which parameters values a Turing bifurcation can occur, given rise to non uniform stationary solutions (patterns) for a reaction-diffusion predator-prey model with variable mortality and Hollyn's type II functional response.
\end{abstract}

\subjclass[2022]{35K57, 92D25}

\keywords{Turing instability, Turing bifurcation, pattern formation, reaction-diffusion systems, predator-prey model}
\maketitle

\section{Introduction}\label{Intro}

Pattern formation in chemical, physical and biological models described by reaction-diffusion systems has been a phenomenon widely investigated and discussed in a large number of scientific publications due to the interest of researchers in the areas of applied mathematics, mathematical biology, engineering chemistry, among others (see for instance \cite{Hu et all, Zhang et al, Jeong et al} and the references therein).
The first published work in this area was carried out by A. M. Turing in 1952 \cite{Turing}, who 
emphasize the role of nonequilibrium diffusion-reaction processes and patterns in biomorphogenesis. Since then, dissipative non-equilibrium mechanisms of spontaneous spatial and spatiotemporal pattern formation in a uniform environment have been of uninterrupted interest in experimental and theoretical biology and ecology (see \cite[Chapter 8]{Malshow Petrovskii and Venturino} and the references therein).
Turing's idea was simple, though not trivial: A stationary state that is locally asymptotically stable in a non-spatial system can become unstable in the corresponding diffusive system and cause stationary solutions (patterns) to appear. This emergent phenomenon of spatial symmetry breaking is known as Turing Instability or Diffusion-Driven Instability (see \cite[Chapter 2]{Grindrod}, \cite{Satnoianu et al 2}) and can be applied to explain pattern formation in the nature, for instance, the appearance of spots on the skin of animals \cite{Jeong et al, Baker et al, Barrios et al}. So, our main goal is to study the Diffusion-Driven Instability phenomenon in a unified way.

In this paper, we investigate the existence of patterns for a reaction-diffusion system with zero-flux (Neumann) boundary conditions (which are usually in pattern formation problems) using the Turing approach. 
In order to state the problems under consideration and the main results of this work, we introduce the general theory to study the effect of diffusion on the stability of the equilibria in a reaction-diffusion model and describe briefly the definitions and concepts related  to  Turing's mechanism to find patterns.\\

\noindent
\textbf{\emph{Turing's mechanism to find patterns:}}
We consider  the following general  reaction-diffusion coupled system of $k$ ($k\geq 1,\; k$ an integer) components (species, concentration of chemicals, etc) which interact in a nonlinear manner and diffuse in $n-$space dimensions with zero-flux across the boundaries
:
\begin{equation}\label{Gensyst1}
	\left \{
	\begin{array}{llll}
		\displaystyle{\frac{\partial w(x,t)}{\partial t}} = D \displaystyle{\Delta w(x,t)}+ G(w(x,t)),\;\;\;x\in\Omega,\;\;t>0,\\
		\displaystyle{\frac{\partial w}{\partial \eta}}(x,t)=0,\;\;\; x \in \partial \Omega,\;\;t>0,	
	\end{array}
	\right.
\end{equation}
where $D=\text{diag}[d_{1},d_{2},\cdots, d_{k}]$ is the diagonal matrix of diffusivities,  with $d_{i}>0,$  $i=1,...,k,$ and $\Delta$ denotes the Laplacian operator on a  bounded, open and connected domain $\Omega \subset \mathbb{R}^{n}$  with  boundary, $\partial \Omega.$  The vector field $\eta(x)$ is the outer unit normal to  $\partial \Omega$ at $x \in \partial \Omega$ and $\displaystyle{\frac{\partial}{\partial \eta}}$ denotes the differentiation in the direction of the outward normal to $\Omega,$ i.e., $\displaystyle{\frac{\partial}{\partial \eta}= \eta \cdot \nabla}$ where $\nabla$ denotes the gradient operator on $\Omega.$ The existence and uniqueness of solutions  for the parabolic system  \eqref{Gensyst1} in the space of continuous functions can be shown (see subsection \ref{eu}) if we assume that the initial data   	
$w(x,0)=\phi(x),\;\;x\in\Omega,$ with $\phi=\left(\phi_{1}, \phi_{2},\cdots, \phi_{k}\right),$  belongs to the set
\begin{equation}\label{Gensyst2}
	X_{\Lambda_{k}}:=\left\{ \phi \in X: \phi(x) \in \Lambda_{k}, x \in \overline{\Omega} \right\},
\end{equation}
where  $X$ is the Banach space 
 \begin{equation}\label{Gensyst3}
 X:=\prod_{i=1}^{k}X_{i},\;\text{with}\;X_{i}=C(\overline{\Omega}, \mathbb{R}),\;i=1,...,k,
 \end{equation}
endowed with the norm $$\|\phi\|:=\displaystyle{\sum_{i=1}^{k}\|\phi_{i}(x)\|_{\infty}},$$ 
and
\begin{equation}\label{Gensystem23}
\Lambda_{k}:=\left\{ v=(v_{1},v_{2},...,v_{k})\in \mathbb{R}^{k}: v_{i}\geq 0,\; i=1,...,k\right\}.
\end{equation}
From the applications point of view,  solutions $w$ of \eqref{Gensyst1} belonging to the set $\Lambda_{k}$ are, in general, relevant.
On the other hand,  in order to have classical solutions  $w=(w_{1}, w_{2},...,w_{k}),$  for system \eqref{Gensyst1}, we  suppose that the nonlinear term $G: \Lambda_{k} \to \mathbb{R}^{k},$ defined by
$$G(w)=\left(G_{1}(w), G_{2}(w), \cdots, G_{k}(w)\right),$$
 satisfies \\
 
 \noindent
 {\bf{H1.}}
 $G\in C^{2}( \mathring{\Lambda}_{k},  \mathbb{R}^{k}),$ and
$G_{i}(w)\geq 0$ whenever $w\in  \Lambda_{k}$ and $w_{i}=0,\;i=1,2,...,k.$\\

 Here $\mathring{\Lambda}_{k}$ denotes the interior of $\Lambda_{k}.$ Finally, we assume that the solutions $w(x,t)$ of \eqref{Gensyst1} are global, i.e. $t\in [0,+\infty).$ 
Next, we summarize the diffusion-driven instability theory, which allow one to analyze the effect of diffusion on the stability of the equilibria  in a reaction-diffusion systems.\\

We said that an equilibrium $ w_{0}=(w_{0}^{1}, w_{0}^{2},\cdots, w_{0}^{k}) \in \Lambda_{k}$ of \eqref{Gensyst1} is \textbf{\textit{Turing (diffusionally) unstable}} if it is an asymptotically stable equilibrium for the kinetic system associated to \eqref{Gensyst1}, i.e.,
\begin{equation}\label{Gensyst4}
	\displaystyle{\frac{dw}{dt}} =  G(w),\;\;\;t>0,
\end{equation}
 but $w_{0}$ is unstable with respect to \eqref{Gensyst1}, see \cite{Satnoianu et al 1} \cite[Chapter 2]{Grindrod} (see also \cite[$\S 10.6$ ]{Okubo}, \cite{Turing, Lizana and Marin} and the references therein).

\begin{rem}
	Note that an equilibrium $w_{0}$ of equation \eqref{Gensyst4}, i.e., $G(w_{0})=0,$ is also a spatial homogeneous equilibrium of equation \eqref{Gensyst1}.
\end{rem}

In \cite[Chapter 2]{Grindrod} the author study the stability  of the homogeneous stationary solution $w_{0}$ of \eqref{Gensyst1} by linearized stability analysis. Specifically, assuming that $w_{0}\in \Lambda_{k}$ is a stable  non-trivial equilibrium of equation \eqref{Gensyst4}, i.e., the Jacobian matrix  $dG(w_{0})$ associated with \eqref{Gensyst4} at $w_{0}$ is stable
$$
	\displaystyle{dG(w_{0}):=\left(
		\begin{array}{cccc}
			\frac{\partial G_{1} (w_{0})}{\partial w_{1}}&	\frac{\partial G_{1} (w_{0})}{\partial w_{2}}& \cdots & 	\frac{\partial G_{1} (w_{0})}{\partial w_{k}} \\
			\frac{\partial G_{2} (w_{0})}{\partial w_{1}}&	\frac{\partial G_{2} (w_{0})}{\partial w_{2}}& \cdots & 	\frac{\partial G_{k} (w_{0})}{\partial w_{k}}\\
			\vdots & \cdots & \ddots & \vdots\\
			\frac{\partial G_{k} (w_{0})}{\partial w_{1}}&	\frac{\partial G_{k} (w_{0})}{\partial w_{2}}& \cdots & 	\frac{\partial G_{k} (w_{0})}{\partial w_{k}}\\	
		\end{array}
		\right),}
$$
one can set $z=w-w_{0},$ and the linearized system \eqref{Gensyst1} around $w_{0}$ is given by
\begin{equation}\label{Gensyst6}
	\left \{
	\begin{array}{llll}
		\displaystyle{\frac{\partial z(x,t)}{\partial t}} = D \displaystyle{\Delta z(x,t)}+ dG(w_{0})z(x,t),\;\;\;x\in\Omega,\;\;t>0,\\
		\displaystyle{\frac{\partial z}{\partial \eta}}(x,t)=0,\;\;\; x \in \partial \Omega,\;\;t>0,
	\end{array}
	\right.
\end{equation}
Therefore, if we denote by $\phi_{j}$ the $j-$th eigenfunction of $-\Delta$ on $\Omega$ with no flux-boundary conditions, that is 
\begin{equation}\label{Gensyst7}
	\displaystyle{\Delta \phi_{j}}+ \lambda_{j}\phi_{j}=0,\; x\in\Omega;\;\;\qquad
	\eta \cdot \nabla \phi_{j}=0,\; x \in \partial \Omega,
\end{equation}
for scalars $\lambda_{j}$ satisfying $0=\lambda_{0}<\lambda_{1} \leq \lambda_{2}\leq \cdots,$
then one can use the  orthogonal basis of $L^{2}(\Omega)$ (see \cite[\& 8.5, pag. 229]{Jurgen}), formed by the set $\{\phi_{j}\}_{j \in \mathbb{N}},$ to solve \eqref{Gensyst6} by expanding the solution $z$   as the Fourier series
\begin{equation}\label{Gensystem8}
	z(x,t)=\sum_{j=0}^{\infty} s_{j}(t) \phi_{j}(x)
\end{equation}
where $s_{j}(t)\in \mathbb{R}^{k},$ for all $j\in \mathbb{N}.$ Hence, substituting \eqref{Gensystem8} in \eqref{Gensyst6}
and equating the coefficients of each eigenfunction $\phi_{j},$ we have that for each $j\in \mathbb{N},$ $s_{j}$ satisfies the ODE
$$\frac{ds_{j}}{dt}=B_{j} s_{j}(t),\quad t>0,$$
where 
\begin{equation}\label{Gensystem9}
	B_{j}=dG(w_{0})-\lambda_{j}D.
\end{equation}
Therefore, the trivial solution $z=0$ ($w=w_{0}$) of \eqref{Gensyst6}  is asymptotically stable if and only if every solution of \eqref{Gensyst6} decays to zero as $t \to +\infty.$ This is equivalent to prove that each $s_{j}(t)$ decays to zero as $t \to +\infty,$ i.e., each matrix $B_{j}$ has $k$ eigenvalues with negative real parts for all $j\in \mathbb{N}.$ On the other hand, if any matrix $B_{j}$ has an eigenvalue with positive real part, then $|s_{j}|$ can grow exponentially and hence so will $z.$ Clearly, in this case, $z=0$ is unstable to arbitrary perturbations which are not orthogonal to the critical eigenmode $\phi_{j}.$ 

\begin{rem}\label{Gensystem13}
	The determination of the pairs $(\phi_{j},\lambda_{j})$ mentioned in \eqref{Gensyst7} is a standard problem of functional analysis, see 
	\cite[pag. 205-208]{Folland}, \cite[\& 8.5, pag. 229]{Jurgen} and \cite{Rodney Biezuner}. We recall that for bounded open domains $\Omega\subset \mathbb{R}^{n},$ the Laplacian operator has a countable infinite number of eigenvalues $\lambda_{j},$ $j\in \mathbb{N},$ each one with finite multiplicity and $\lambda_{j} \to +\infty$ as $j\to +\infty$ (see particular examples in subsection \ref{patt}). Furthermore, the differential operator $-\Delta,$ with no-flux boundary conditions, is self-adjoint in $L^{2}(\Omega),$ and for each $j\geq 1$ it is easy to show that 
	$$\lambda_{j}=\displaystyle{\frac{\int_{\Omega} |\nabla \phi_{j}|^{2}dx}{\int_{\Omega} \phi_{j}^{2}dx}>0.}$$
\end{rem}
If parameters $d_{1}, d_{2}, ..., d_{k}$ are such that some $B_{j}$ has an eigenvalue with zero real part, then the stability of the trivial solution $z=0$  of system \eqref{Gensyst6}  will switch and this shall reflect (under appropriate conditions) a bifurcation of some inhomogeneous equilibrium from the equilibrium $w\equiv w_{0}$ for \eqref{Gensyst1}. Since, in general, a diffusive coefficient is considered  as a bifurcation parameter, we  said that an equilibrium $ w_{0} \in \Lambda_{k}$ of system \eqref{Gensyst1} undergoes a \textbf{\textit{Turing bifurcation}} at $d^{\ast}\in (0, +\infty)$ if the solution $w_{0}$ change its stability at $d^{\ast}$ and in some neighborhood of $d^{\ast}$ there exists a one-parameter family of non-constant stationary solution of system \eqref{Gensyst1}.\\

\begin{rem}
	Recently, the authors in \cite{Krause et al} have shown that there exist reaction-diffusion models satisfying the Turing instability conditions but robustly exhibit only transient patterns. This point out the importance of finding sufficient conditions for the existence of a patterned state in reaction-diffusion models.
\end{rem}

\noindent
\textbf{\emph{Problems under consideration:}}
From \eqref{Gensystem9} and the analysis carried out above  the following questions arise:

\begin{itemize}
	\item[$1)$] Which properties must the Jacobian matrix  $dG(w_{0})$ satisfy  to ensure  the equilibrium $ w_{0} \in \Lambda_{k}$ of \eqref{Gensyst1} is Turing  unstable?.
	
	\item[$2)$] How one must select the pair $(\phi_{j},\lambda_{j})$ to allow the matrix $B_{j}$ given in \eqref{Gensystem9} has an eigenvalue with positive real part? What properties should have the eigenvalue $\lambda_{j}$ to assure that the equilibrium $ w_{0} \in \Lambda_{k}$ of \eqref{Gensyst1} is Turing  unstable?	
	
	\item [$3)$] How exactly should the diffusive parameters $d_{1}, d_{2}, ..., d_{k}$ be chosen to ensure that equilibrium $ w_{0} \in \Lambda_{k}$ of \eqref{Gensyst1} is Turing  unstable and, furthermore, undergoes a Turing bifurcation?
\end{itemize}

In the literature, there is a wide range of references answering the first question since   Turing instability in reaction-diffusion systems can be recast in terms of matrix stability (see \eqref{Gensystem9}), we refer the reader to \cite{Cross, Hershkowitz, Satnoianu et al 1, Satnoianu et al 2, Hadeler and Ruan} and the references therein to delve deeper into this topic.
In the Appendix A, we review some important results on stability of real matrices. Special emphasis is given to the concept of excitable matrix, which will be  fundamental to prove the main results of this work. Necessary and sufficient conditions to determine when a $2\times2$ or $3\times 3$ matrix is excitable have been established in \cite{Cross} (see Theorems \ref{strsta2x2} and \ref{strsta3x3}) obtaining, in this way, four possibilities of sign structures for the Jacobian matrix $dG(w_{0})$ in the $2\times 2$ case (see Remark \ref{excita2x21}). There are only sufficient conditions for real $k\times k-$matrices, with $k\geq 4,$ to be excitable, see Remark \ref{necexcitable}. However, as far as we know, there are not general results regarding  questions $2)$ and $3).$ \\

\noindent
\textbf{\emph{Main results:}} 
In the present manuscript, we apply a bifurcation from a  simple eigenvalue result for operators in Banach spaces (see \cite[Theorem 13.5]{Smoller}), to find a practical criterion, given specific conditions on the diffusive parameters, to show that the non-trivial equilibrium $w_{0}$ of \eqref{Gensyst1} undergoes a Turing bifurcation and hence to obtain the existence of patterns in reaction-diffusion models (see \eqref{Gensyst1}) of two components ($k=2$). Therefore, we were able to extend the techniques used by the authors in \cite{Lizana and Marin, Cavani and Farkas II}
to a wide class of two-components reaction-diffusion models. For this, we assume that 
\begin{equation}\label{nececond}
w_{0}=(w_{0}^{1}, w_{0}^{2}) \in\mathring{\Lambda}_{2}
\end{equation}
is a nontrivial equilibrium of equation \eqref{Gensyst4} with $k=2,$ i.e., $w_{0}^{i} \neq 0$ for $i=1,2,$ and set
\begin{equation}\label{Gensystem10}
	\displaystyle{A:=\left(
		\begin{array}{cccc}
			a_{11}&a_{12} \\
			a_{21}&a_{22}\\
		\end{array}
		\right)
		=\left(
		\begin{array}{cccc}
			\frac{\partial G_{1} (w_{0})}{\partial w_{1}}&	\frac{\partial G_{1} (w_{0})}{\partial w_{2}} \\
			\frac{\partial G_{2} (w_{0})}{\partial w_{1}}&	\frac{\partial G_{2} (w_{0})}{\partial w_{2}}\\
		\end{array}
		\right)
		=dG(w_{0}).}
\end{equation}
Also, $D=\text{diag}[d_{1},d_{2}]$ with $d_{i}>0,$  $i=1,2,$ and  for $\lambda \in \mathbb{R}_{+}$ given, we define the Hyperbola 	$H_{\lambda}$ in the $d_{1}d_{2}-$plane by
\begin{equation}\label{Gensystem11}
	H_{\lambda}:=\{(d_{1},d_{2})\in \mathbb{R}^{2} :  \text{det}(A)-a_{11}\lambda d_{2} - a_{22}\lambda d_{1}+\lambda^{2}d_{1}d_{2}= 0\}.
\end{equation}

\begin{rem}
Condition \eqref{nececond}	is necessary to prove Theorems \ref{Gensystem12}-\ref{Gensystem18}. Note that authors in \cite{Krause et al} consider  an equilibrium $w_{0}$ in the border $\partial \Lambda_{2}.$ On the other hand, they used periodic boundary conditions instead of Neumann boundary conditions. 
\end{rem}

The following results state that if a matrix $A$ is excitable, and it is known a simple eigenvalue  $\lambda_{l},$ (so, its corresponding eigenfunction $\phi_{l}$) of the Laplacian operator for some $l\in \mathbb{N}^{\ast},$ then we can exchange the stability of the matrix $B_{l}$ (see \eqref{Gensystem9} with $k=2$) provided that the diffusive coefficients are chosen appropriately on the $d_{1}d_{2}-$plane.  
The choice of the diffusive bifurcation parameter shall depend on the conditions for  the matrix A to be excitable (see \eqref{excita2x2}).
Our first main result consider $d_{1}$ as bifurcation parameter and reads as follows.

\begin{thm}[Criterion I]\label{Gensystem12}
	Assume that system \eqref{Gensyst1} satisfies {\bf{H1}} (with $k=2$) and $A\in \mathcal{M}_{2\times 2}(\mathbb{R})$  given in \eqref{Gensystem10}	satisfies \eqref{strsta2x21}, \eqref{strsta2x22} and $a_{22}>0.$ Let $\lambda_{j}$ be the eigenvalues,  with respective eigenfunctions $\phi_{j},$ of operator $-\Delta$ on $\Omega$ with no flux-boundary conditions for $j=0,1,2,...$ (see \eqref{Gensyst7} and Remark \eqref{Gensystem13}).
	Suppose that for some $l\in \mathbb{N},$ with $l \geq 2,$ the eigenvalue $\lambda_{l}$ is simple and consider the intersection points   $(d_{1}',d_{2}')\in H_{\lambda_{l+1}} \cap H_{\lambda_{l} },$ $(d_{1}'',d_{2}'')\in H_{\lambda_{l}} \cap H_{\lambda_{l-1} }$ (see \eqref{Gensystem11} and Figure \ref{Bifurca22}) that are located in the first quadrant of the $d_{1}d_{2}-$plane (Note that $0<d_{2}'$ and  $d_{2}''  <\frac{a_{22}}{\lambda_{l}}).$ 	If
	\begin{equation}\label{Gensystem14}
		d_{2}'<d_{2}<d_{2}'' 
	\end{equation}	
	then,  at $\displaystyle{d_{1}^{\ast}=\frac{a_{11}\lambda_{l}d_{2} -\text{det}(A)}{\lambda_{l}(\lambda_{l}d_{2}-a_{22})}},$
	the uniform steady-state solution $w_{0}\in \mathring{\Lambda}_{2}$ of 
	\eqref{Gensyst1} with $k=2,$ undergoes a Turing bifurcation. Furthermore,
	\begin{equation}\label{Gensystem15}
w(x,s)=w_{0}+s\cdot v_{1l}\phi_{l}(x)+O(s^{2})
	\end{equation}
	are non-uniform stationary solutions of \eqref{Gensyst1} with $k=2,$   $d_{1}\equiv d_{1}(s),$ $s\in (-\zeta, \zeta),$ for some positive $\zeta$ small enough 
	and  $v_{1l}$ (see \ref{eigenvectorsform1} with $j=l$ and $d_{1}=d_{1}^{\ast}$) is the eigenvector associated with eigenvalue $\lambda_{1l}=0$ of the matrix
	\begin{equation}\label{Gensystem16}
		B_{l}=A-\lambda_{l}\left(
		\begin{array}{cc}
			d_{1}^{\ast} & 0 \\
			0& d_{2}\\	
		\end{array}
		\right).
	\end{equation}
\end{thm}

\begin{rem}\label{Gensystem17}
	Theorem \ref{Gensystem12} gives the exact intervals where diffusive parameters should variate in order to exchange the stability of only one matrix, $B_{l},$ and $\lambda_{l}$ is the unique eigenvalue allowing this property, i.e., $B_{j},$ for $j=0,1,2\cdots, l-1,l+1,...,$ remains as stable matrices. Furthermore, it is possible to give explicitly the interval where the diffusion coefficient $d_{2}$ should be taken, i.e.,
$$d_{2}'=\frac{\text{det}(A)(\lambda_{l}+\lambda_{l+1})-\sqrt{[\text{det}(A)(\lambda_{l}+\lambda_{l+1})]^{2}-4a_{11}a_{22}\text{det}(A)\lambda_{l}\lambda_{l+1}}}{2a_{11}\lambda_{l}\lambda_{l+1}},$$

$$d_{2}''=\frac{\text{det}(A)(\lambda_{l}+\lambda_{l-1})-\sqrt{[\text{det}(A)(\lambda_{l}+\lambda_{l-1})]^{2}-4a_{11}a_{22}\text{det}(A)\lambda_{l}\lambda_{l-1}}}{2a_{11}\lambda_{l}\lambda_{l-1}}.$$
Also, it is easy to show that under hypotheses of Theorem \ref{Gensystem12}, for fixed $d_{2}$  satisfying \eqref{Gensystem14}, a neighborhood of $d_{1}^{\ast}$ that preserves the desired structure for eigenmodes is an  open interval with upper bound $ d_{1}''',$ where $$d_{1}''':=\min\left\{\frac{a_{11}\lambda_{l+1}d_{2}-\text{det}(A)}{\lambda_{l+1}(\lambda_{l+1} d_{2} -a_{22})},\; \frac{a_{11}\lambda_{l-1}d_{2}-\text{det}(A)}{\lambda_{l-1}(\lambda_{l-1}d_{2} -a_{22})} \right\}.$$ 
	Thus, for numerical simulations,  the neighborhood of $d_{1}^{\ast}$ can be considered as a sub-interval of $(d_{1}',d_{1}''').$ 
	If $\lambda_{1}$ is simple we can also exchange the stability of $B_{1}.$ For this, we fix $d_{2}$ satisfying $\frac{a_{22}}{\lambda_{2}}\leq d_{2}<\frac{a_{22}}{\lambda_{1}}.$ In this way, at $\displaystyle{d_{1}^{\ast}=\frac{a_{11}\lambda_{1}d_{2} -\text{det}(A)}{\lambda_{1}(\lambda_{1}d_{2}-a_{22})}},$ the uniform steady-state solution $w_{0}\in \mathring{\Lambda}_{2}$ of 
	\eqref{Gensyst1} with $k=2,$ undergoes a Turing bifurcation. 
\end{rem}

The following result study the stability of the bifurcating solution \eqref{Gensystem15}.

\begin{prop}\label{stability}
 Assume that the conditions of Theorem \ref{Gensystem12} are satisfied and	let $(d_{1}(s), w(x,s))$ be the one-parameter family of bifurcating solutions given by the formula \eqref{Gensystem15}. Suppose that  $d_{1}'(0)\neq 0,$ and the eigenvalues, say $\gamma(s),$ of the nonhomogeneous steady state bifurcating from the critical value $\lambda_{1l}=0$ are non-zero for small $|s|\neq 0.$ Then, if $d_{1}(s)<d_{1}^{\ast}$ the corresponding solution $w(x,s)$ is unstable and if $d_{1}(s)> d_{1}^{\ast},$ the corresponding solution $w(x,s)$ is stable.
\end{prop}

Next we shall prove that if the excitable matrix  $A$ given in \eqref{Gensystem10}	satisfies $a_{11}>0,$ then one can consider $d_{2}$ as a bifurcation parameter to determine nonhomogeneous stationary solutions of \eqref{Gensyst1} with $k=2.$ Our second main result regarding Turing bifurcation reads as follows. 

\begin{thm}[Criterion II]\label{Gensystem18}
	Assume that system \eqref{Gensyst1} satisfies {\bf{H1}} (with $k=2$) and  $A\in \mathcal{M}_{2\times 2}(\mathbb{R})$  given in \eqref{Gensystem10} satisfies \eqref{strsta2x21}, \eqref{strsta2x22} and $a_{11}>0.$ 
	Let $\lambda_{j}$ be the eigenvalues,  with respective eigenfunctions $\phi_{j},$ of operator $-\Delta$ on $\Omega$ with no flux-boundary conditions for $j=0,1,2,...$ (see \eqref{Gensyst7} and Remark \eqref{Gensystem13}).
	Suppose that for some $l\in \mathbb{N},$ with $l \geq 2,$ the eigenvalue $\lambda_{l}$ is simple and consider the intersection points   $(d_{1}',d_{2}')\in H_{\lambda_{l+1}} \cap H_{\lambda_{l} },$ $(d_{1}'',d_{2}'')\in H_{\lambda_{l}} \cap H_{\lambda_{l-1} }$ (see \eqref{Gensystem11} and Figure \ref{Bifur2}) that are located in the first quadrant of the $d_{1}d_{2}-$plane  (Note that  $0<d_{1}'$ and $d_{1}''  <\frac{a_{11}}{\lambda_{l}}$).
	If 
	\begin{equation}\label{Gensystem19}
		d_{1}'<d_{1}<d_{1}''
	\end{equation}
	then, at
	$\displaystyle{d_{2}^{\ast}=\frac{a_{22}\lambda_{l}d_{1} -\text{det}(A)}{\lambda_{l}(\lambda_{l}d_{1}-a_{11})}}$ the uniform steady-state solution $w_{0}\in \mathring{\Lambda}_{2}$ of 
	\eqref{Gensyst1} with $k=2,$ undergoes a Turing bifurcation. Furthermore,
	\begin{equation}\label{Gensystem20}
		w(x,s)=w_{0}+s\cdot v_{1l}\phi_{l}(x)+O(s^{2})
	\end{equation}
	are non-uniform stationary solutions of \eqref{Gensyst1} with $k=2,$   $d_{2}\equiv d_{2}(s),$ $s\in (-\zeta, \zeta),$ for some positive $\zeta$ small enough 
	and  $v_{1l}$ (see \ref{eigenvectorsform2} with $j=l$ and $d_{2}=d_{2}^{\ast}$) is the eigenvector associated with eigenvalue $\lambda_{1l}=0$ of the matrix 
$$
		B_{l}=A-\lambda_{l}\left(
		\begin{array}{cc}
			d_{1} & 0 \\
			0& d_{2}^{\ast}\\	
		\end{array}
		\right).
$$
\end{thm}

\begin{rem}\label{crtII}
One can compute explicitly the interval \eqref{Gensystem19} where the diffusive coefficient $d_{1}$ should be taken, i.e.,
$$d_{1}'=\frac{\text{det}(A)(\lambda_{l}+\lambda_{l+1})-\sqrt{[\text{det}(A)(\lambda_{l}+\lambda_{l+1})]^{2}-4a_{22}a_{11}\text{det}(A)\lambda_{l}\lambda_{l+1}}}{2a_{22}\lambda_{l}\lambda_{l+1}},$$

$$d_{1}''=\frac{\text{det}(A)(\lambda_{l}+\lambda_{l-1})-\sqrt{[\text{det}(A)(\lambda_{l}+\lambda_{l-1})]^{2}-4a_{22}a_{11}\text{det}(A)\lambda_{l}\lambda_{l-1}}}{2a_{22}\lambda_{l}\lambda_{l-1}}.$$
For numerical simulations, it is easy to show that under hypotheses of Theorem \eqref{Gensystem18}, for  fixed $d_{1}$ satisfying condition \eqref{Gensystem19}, a neighborhood of $d_{2}^{\ast}$ that preserves the desired structure for the eigenmodes is an open interval with upper bound $d_{2}''',$  where 
$$
		d_{2}''':=\min\left\{\frac{a_{22}\lambda_{l+1}d_{1}-\text{det}(A)}{\lambda_{l+1}(\lambda_{l+1} d_{1} -a_{11})},\; \frac{a_{22}\lambda_{l-1}d_{1}-\text{det}(A)}{\lambda_{l-1}(\lambda_{l-1}d_{1} -a_{11})} \right\}.
$$

On the other hand,  if $\lambda_{1}$ is a simple eigenvalue then, one can exchange the stability of $B_{1}$ by requesting $\frac{a_{11}}{\lambda_{2}}\leq d_{1}<\frac{a_{11}}{\lambda_{1}}.$ Therefore,  at $\displaystyle{d_{2}^{\ast}=\frac{a_{22}\lambda_{1}d_{1} -\text{det}(A)}{\lambda_{1}(\lambda_{1}d_{1}-a_{11})}},$ the uniform steady-state solution $w_{0}\in \mathring{\Lambda}_{2}$ of 
\eqref{Gensyst1} with $k=2,$ undergoes a Turing bifurcation. 
\end{rem}

Next, we study the stability of the bifurcating solution \eqref{Gensystem20}.

\begin{prop}\label{stability2}
Assume that the conditions of Theorem \ref{Gensystem18} are satisfied and	let $(d_{2}(s), w(x,s))$ be the one-parameter family of bifurcating solutions given by the formula \eqref{Gensystem20}.  Suppose that $d_{2}'(0)\neq 0,$ and that the eigenvalues, say $\beta(s),$ of the nonhomogeneous steady state bifurcating from the critical value $\lambda_{1l}=0$ are non-zero for small $|s|\neq 0.$ Then, if $d_{2}(s)<d_{2}^{\ast}$ the corresponding solution $w(x,s)$ is unstable and if $d_{2}(s)> d_{2}^{\ast},$ the corresponding solution $w(x,s)$ is stable.
\end{prop}

The interaction of at least two species with different diffusion coefficients can give rise to this kind of spatial structure,  see  \cite[Part III, Chapter 9]{Malshow Petrovskii and Venturino}.
In this paper, we apply Criterion II (Theorem \eqref{Gensystem18}), to study the existence of patterns for a reaction-diffusion  predator-prey model with variable mortality and Hollyn's type II functional response (see Section \ref{app}). It is known that patterns produced by the Turing mechanism can be sensitive to domain shape \cite{Bunow et al} and, therefore, it is important to investigate the robustness of the above patterns to changes in geometry. In particular, we investigate how the pattern varies due to changes in domain shape and dimension through numerical simulations.\\

The paper is organized as follows. In Section \ref{preliminary}, we prove the well-posedness of system \ref{Gensyst1} (see Subsection \ref{eu}) and we study the behavior of the eigenvalues associated to the matrix $B_{j}$ given in \eqref{Gensystem9} when $k=2$ (see Subsection \ref{eigenbeha}). We show the Theorems \ref{Gensystem12} and  \ref{Gensystem18} in Section \ref{Turing Bufurcation}.  Section \ref{stab}  is dedicated to study the stability of the non-uniform stationary solutions \eqref{Gensystem15} and \eqref{Gensystem20} of system \eqref{Gensyst1} with $k=2$ that arise from the bifurcation of the homogeneous steady state $w_{0}$ given in \eqref{Gensystem10}. Finally, in Section \ref{app},  we apply our results and perform some numeric simulations to prove the existence of patterns for a reaction-diffusion predator-prey model with variable mortality and Hollyn's type II functional response.

\section{Preliminary Results}\label{preliminary}
In this section, we will recall general results regarding the existence and uniqueness of solutions of reaction-diffusion systems. Such results will be used in section \ref{wellpose} to show that   \eqref{Gensystem} generates a dynamical system which is  well-posed on the Banach space $X$ given in \eqref{Gensyst3}.  

\subsection{Well-posedness for Parabolic systems with Newman boundary conditions}\label{eu}
  For the sake of completeness, we summarize the following general results regarding the existence and uniqueness of solutions for parabolic systems with Neumann boundary conditions in the space of continuous functions \cite[Chapter 7]{Hal Smith} (see also \cite{Lizana and Marin}). For this, $\Omega$ will denote a bounded, open and connected subset of  $\mathbb{R}^{n}$ with a piecewise smooth boundary, $\partial \Omega.$ Let $C(\overline{\Omega})=
C(\overline{\Omega},\mathbb{R})$ be the Banach space of continuous functions in $\overline{\Omega}$ endowed with the usual supremum norm denoted by $\|\cdot\|_{\infty}.$ Let $\Delta, \;\nabla,\; D,$ and $\displaystyle{\frac{\partial}{\partial \eta}}$ be defined as in \eqref{Gensyst1} and
consider the general reaction-diffusion system
\begin{equation}\label{Gensystem1}
	\left \{
	\begin{array}{llll}
		\displaystyle{\frac{\partial w(x,t)}{\partial t}} = D \displaystyle{\Delta w(x,t)}+ G(x,w(x,t)),\;\;\;x\in\Omega,\;\;t>0,\\
		\displaystyle{\frac{\partial w}{\partial \eta}}(x,t)=0,\;\;\; x \in \partial \Omega,\;\;t>0,		\\
		w(x,0)=\phi(x),\;\;x\in\Omega,
	\end{array}
	\right.
\end{equation}
where $w=(w_{1}, w_{2}, \cdots , w_{k}),$ $\phi=\left(\phi_{1}, \phi_{2},\cdots, \phi_{k}\right),$ and the nonlinear term $G:\overline{\Omega} \times \Lambda_{k} \to \mathbb{R}^{k}$ is defined by
 $$G(x,w)=\left(G_{1}(x,w), G_{2}(x,w), \cdots, G_{k}(x,w)\right).$$
 Let $A_{i}^{0}$ be the differential operator 
$$A_{i}^{0}w_{i}=d_{i}\Delta w_{i},$$ defined on the domain $D(A_{i}^{0}) \subset X_{i}$ (see \eqref{Gensyst3}) given by
$$D(A_{i}^{0})=\left\{w_{i}\in C^{2}(\Omega) \cap C^{1}(\overline{\Omega}): A_{i}^{0} w_{i} \in C^{1}(\overline{\Omega}), 	\displaystyle{\frac{\partial w_{i}(x)}{\partial \eta}}=0,\;\;\; x \in \partial \Omega \right\}.$$

The closure $A_{i}$ of $A_{i}^{0}$ in $X_{i}$ generates an analytic semigroup of bounded linear operators $T_{i}(t)$ for $t\geq 0$ (see \cite[Theorem 2.4]{Mora}) such that $w_{i}(t)=T_{i}(t)\phi_{i}$ is the solution of the abstract linear differential equation in $X_{i}$ given by
\begin{equation}\label{LLSystem}
	\left \{
	\begin{array}{llll}
		w'_{i}(t)=A_{i}w_{i}(t),\;\; t>0,\\
		w_{i}(0)=\phi_{i}\in D(A_{i}).
	\end{array}
	\right.
\end{equation}
An additional property of the semigroup is that for each $t>0,$ $T_{i}(t)$ is a compact operator. In the language of partial differential equations $$w_{i}(x,t)=[T_{i}(t)\phi_{i}](x)$$
is a classical solution of the initial boundary value problem 
$$
	\left \{
	\begin{array}{llll}
		\displaystyle{\frac{\partial w_{i}(x,t)}{\partial t}} = d_{i} \displaystyle{\Delta w_{i}(x,t)},\;\;\;x\in\Omega,\;\;t>0,\\
		\displaystyle{\frac{\partial w_{i}}{\partial \eta}}(x,t)=0,\;\;\; x \in \partial \Omega,\;\;t>0,		\\
		w_{i}(x,0)=\phi_{i}(x),\;\;x\in\Omega.
	\end{array}
	\right.
$$
Therefore, $T(t):X\to X$ defined by $\displaystyle{T(t):=\prod_{i=1}^{k}T_{i}(t)}$ is a semigroup of operators on $X$ (see \eqref{Gensyst3}) generated by the operator $\displaystyle{A:=\prod_{i=1}^{k} A_{i}},$ defined on $\displaystyle{D(A):=\prod_{i=1}^{k}D(A_{i})},$ and $w(x,t)=[T(t)\phi](x)$ is the solution of the linear system 
$$
	\left \{
	\begin{array}{llll}
		\displaystyle{\frac{\partial w(x,t)}{\partial t}} = D \displaystyle{\Delta w(x,t)},\;\;\;x\in\Omega,\;\;t>0,\\
		\displaystyle{\frac{\partial w}{\partial \eta}}(x,t)=0,\;\;\; x \in \partial \Omega,\;\;t>0,		\\
		w(x,0)=\phi(x),\;\;x\in\Omega.
	\end{array}
	\right.
$$

Note that if we request the nonlinear term $G$ to be twice continuously differentiable in 
$\overline{\Omega} \times \mathring{\Lambda}_{k}$ then, we can define the map $\displaystyle{[f(\phi)](x)=G(x,\phi(x))},$ which maps $X$ into itself, and equation 
\eqref{Gensystem1} can be viewed as the abstract ODE in $X$ (see \eqref{Gensyst3}) given by 

\begin{equation}\label{Gensystem4}
	\left \{
	\begin{array}{llll}
	w'(t)=Aw(t)+f(w(t)),\;\;t>0,\\
		w(0)=\phi.
	\end{array}
	\right.
\end{equation}
While a solution $w(t)$ of \eqref{Gensystem4} can be obtained under the restriction that $\phi \in D(A),$ a so-called mild solution can be obtained for every $\phi \in X$ by requiring only that $w(t)$ is a continuous solution of the following integral equation

$$u(t)=T(t)\phi+\int_{0}^{t}T(t-s)f(u(s))ds,\;\; t\in[0,t_{0}),$$
where $t_{0}=t_{0}(\phi)\leq \infty.$ Restricting our attention to functions $\phi$ in the set $X_{\Lambda_{k}}$ (see \eqref{Gensyst2}), and supposing that 
$G:\overline{\Omega} \times \Lambda_{k} \to \mathbb{R}^{k}$ satisfies 
$G_{i}(x,w)\geq 0$ whenever $(x,w)\in \overline{\Omega} \times \Lambda_{k}$ and $w_{i}=0,$ then Corollary 3.3 in \cite[pp. 129]{Hal Smith} implies the Nagumo condition for the positive invariance of $\Lambda_{k}$ (see \eqref{Gensystem23}), i.e.,
\begin{equation}\label{Nagumo}
\lim_{h\to 0^{+}} h^{-1} \text{dis}(\Lambda_{k}, w+h G(x,w))=0, \;\;(x,w)\in \overline{\Omega} \times \Lambda_{k}.
\end{equation}
Furthermore, a direct application of the strong parabolic maximum principle (see \cite[Theorem 2.1]{Hal Smith}) show that the linear semigroup $T(t)$ leaves $X_{\Lambda_{k}}$ positively invariant, that means,
\begin{equation}\label{PI}
T(t)X_{\Lambda_{k}} \subset X_{\Lambda_{k}}, \;\;t\geq 0.
\end{equation}
From conditions \eqref{Nagumo} and \eqref{PI} we infer the following result.

\begin{thm}\label{GenerWP}
	Assume that the nonlinear term $G$ present in \eqref{Gensystem1} satisfies $G\in C^{2}(\overline{\Omega} \times \mathring{\Lambda}_{k},  \mathbb{R}^{k}),$ 
	and $G_{i}(x,w)\geq 0$ whenever $(x,w)\in \overline{\Omega} \times \Lambda_{k}$ and $w_{i}=0.$ Then for each $\phi \in X_{\Lambda_{k}}$  (see \eqref{Gensyst2}),
	the reaction-diffusion system \eqref{Gensystem1} has an unique noncontinuable mild solution $w(t)=w(t,\phi)\in X_{\Lambda_{k}}$ defined on $[0,\sigma),$ where $\sigma=\sigma(\phi)\leq \infty.$ Furthermore, the following properties hold:
	\begin{itemize}
		\item[i)] $w(t)$ is continuously differentiable on $(0,\sigma),$ $w(t)\in D(A)$ and $w(t)$ satisfies \eqref{Gensystem4} on $(0, \sigma).$
		
		\item[ii)] $w(x,t)=[w(t)](x)$ is a classical solution of \eqref{Gensystem1}.

		\item[iii)] If $\sigma<+\infty,$ then $\|w(t)\|\to \infty$ as $t\to \infty.$
		
		\item[iv)] If $\sigma (\phi)=+\infty$ for all $\phi \in X_{\Lambda_{k}},$ then $\Phi_{t}(\phi)=w(t,\phi)$ is a semiflow on $X_{\Lambda_{k}}.$
	\end{itemize}
\end{thm}
\begin{proof}
See 	\cite[Theorem 3.1]{Hal Smith}
\end{proof}

Note that Theorem \ref{GenerWP} and hypothesis {\bf{H1}} imply that system \eqref{Gensyst1} with initial data $\phi \in X_{\Lambda_{k}}$ is well-posed. 
Next, we analyze the stability of of the matrix $B_{j}$ given in \eqref{Gensystem9} in the particular case $k=2.$

 \subsection{Behavior of eigenvalues of the matrix $B_{j}$ given in \eqref{Gensystem9} when $k=2$.}\label{eigenbeha}
\begin{figure}[h!]
\includegraphics[width=5cm]{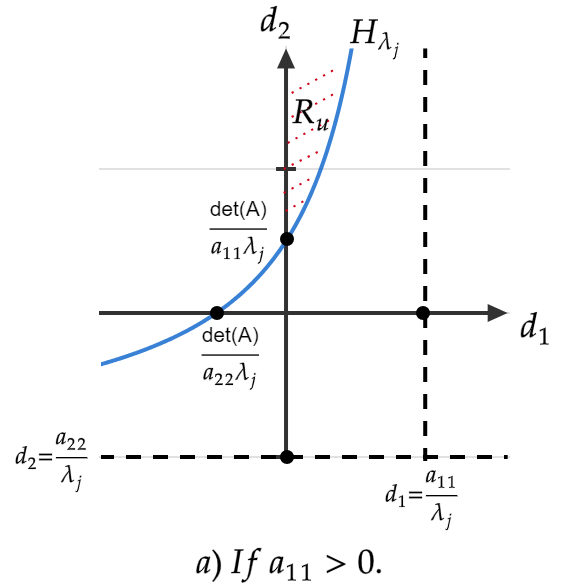}
\includegraphics[width=9cm]{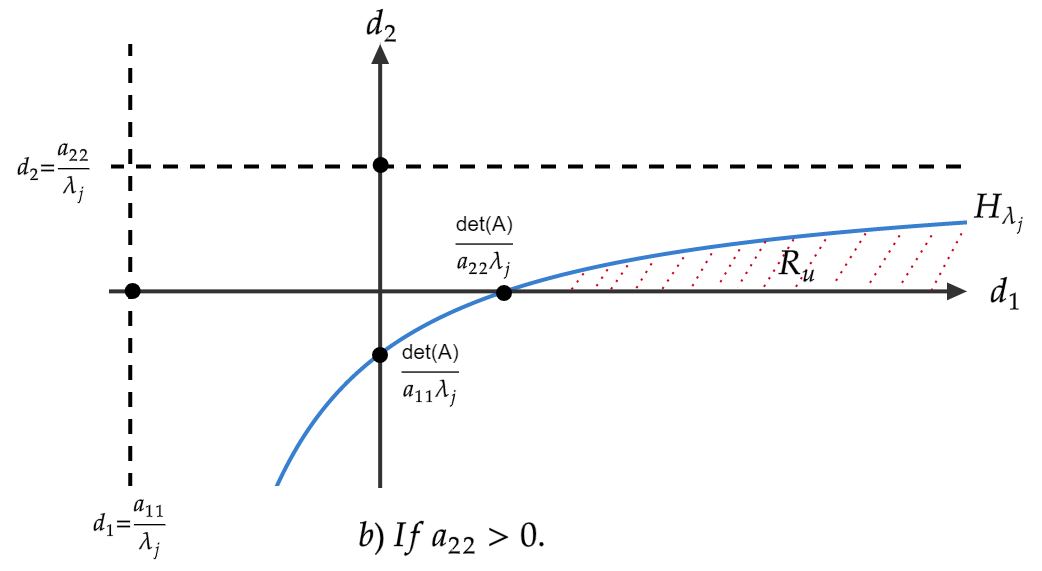}
\caption{Region $R_{u}$ in the $d_{1}d_{2}-$plane.}\label{casoa22}
\end{figure}

\begin{figure}[h!]
 	\includegraphics[width=9cm]{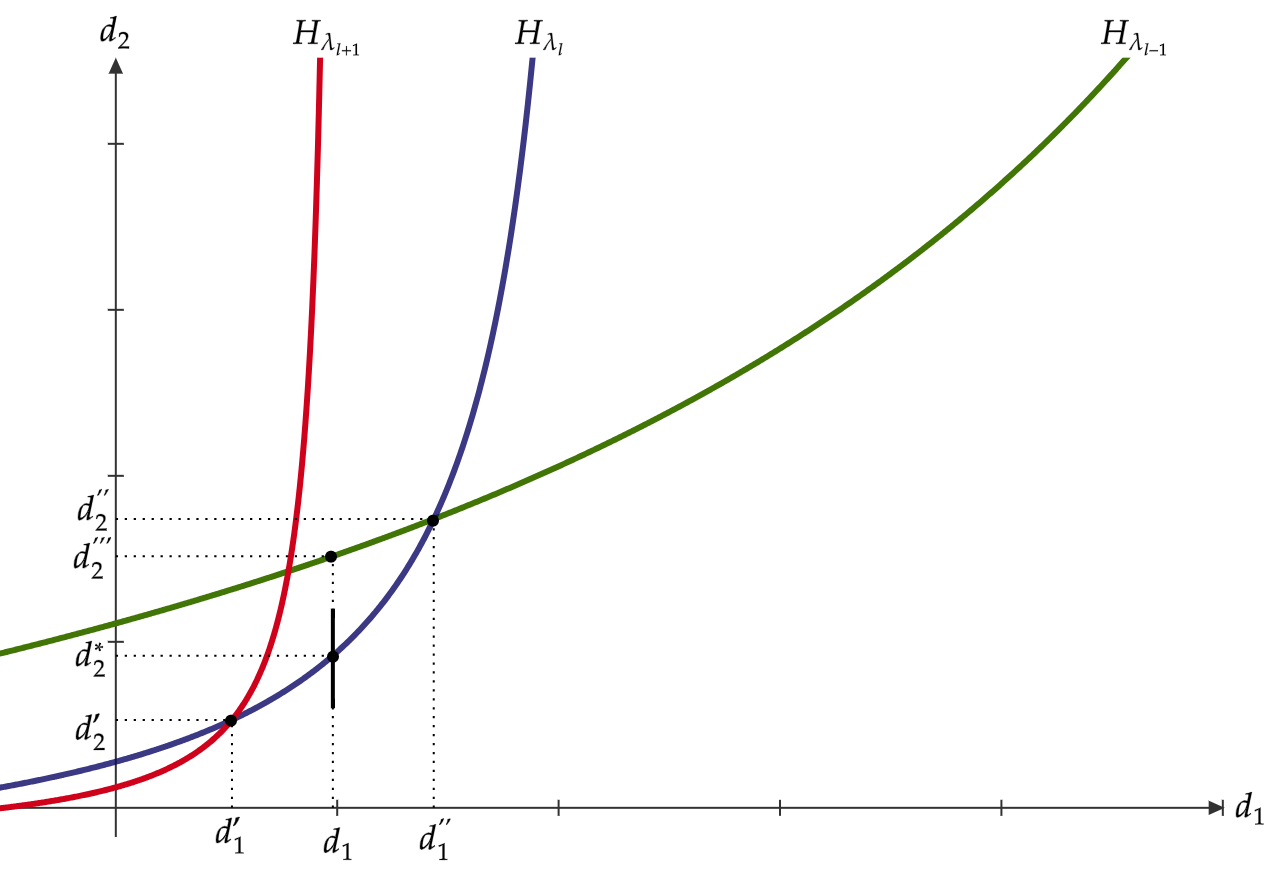}
 	\caption{Turing Bifurcation: when $d_{1}'<d_{1}<d_{1}'',$ the uniform steady state solution $w_{0}\in \mathring{\Lambda}_{2}$ of \eqref{Gensyst1} with $k=2$ undergoes a Turing Bifurcation at $d_{2}=d_{2}^{\ast}.$  Here, $\lambda_{j}$ denote the $j$th eigenvalue of $-\Delta$ on $\Omega$ with no flux-boundary conditions (see \eqref{Gensyst7}) and the hyperbola $H_{\lambda_{j}},$ $j=1,2,...,l-1,l, l+1,...$ is defined as in \eqref{Gensystem11} (see Figure \eqref{casoa22} part $a)$).}\label{Bifur2}
 \end{figure}

Here we fix $j\in \mathbb{N}^{\ast}$ and analyze the behavior of the eigenvalues of the matrix
\begin{equation}\label{ma2}
B_{j}=A-\lambda_{j}D,
\end{equation} 
where $A\in \mathcal{M}_{2\times 2}(\mathbb{R})$ and $D$ are defined as in \eqref{Gensystem10}, and $\lambda_{j}>0$ is given by \eqref{Gensyst7}.  The matrix $A,$ is assumed to be  excitable (see Theorem \ref{strsta2x2} in the Appendix A).  In other words, in this subsection we shall study the stability of nontrivial equilibrium $w_{0}=(w_{0}^{1}, w_{0}^{2}) \in \mathring{\Lambda}_{2}$ with respect to the system \eqref{Gensyst1} with $k=2.$
The eigenvalues of the matrix $B_{j}$ are given by the roots of the characteristic polynomial 
\begin{equation}\label{charpolik=2}
	P_{B_{j}}(\rho)=
	\rho^{2}-\text{Tr}(B_{j}) \rho +\text{det}(B_{j})
\end{equation}
From Theorem \ref{strsta2x2}, we infer $-\text{Tr}(B_{j})=-\left[\text{Tr}(A)-\lambda_{j}\text{Tr}(D)\right]>0.$ The Routh-Hurwitz criterion imply that, for the Turing instability to occur, it should be satisfied that $\text{det}(B_{j})\leq 0,$  where $\text{det}(B_{j})=(a_{11}-\lambda_{j}d_{1}) (a_{22}-\lambda_{j}d_{2})-a_{12} a_{21}.$
Based in the two possibilities expressed in \eqref{excita2x2} and Remark \ref{excita2x21}, we can depict the plane region $R_{u},$  where we must take the positive diffusive coefficients $(d_{1},d_{2})$ such that homogeneous stationary solution  $w_{0}$ is unstable to the system \eqref{Gensyst1} with $k=2.$ Analytically, the region $R_{u}$ is described by
$$R_{u}=\{(d_{1},d_{2})\in \mathbb{R}_{+}^{2} :  \text{det}(A)-a_{11}\lambda_{j} d_{2} - a_{22}\lambda_{j} d_{1}+\lambda_{j}^{2}d_{1}d_{2}< 0\},$$
and for such $\lambda_{j}>0,$  we consider the Hyperbola  $H_{\lambda_{j}}$ in the $d_{1}d_{2}-$plane defined as in \eqref{Gensystem11} (see Figure \ref{casoa22}).

The following result states that if the positive diffusive coefficients cross the Hyperbola $H_{\lambda_{j}}$ from the instability region $R_{u}$  then, by continuity, a real eigenvalue of the matrix $B_{j}$ given in \eqref{ma2} cross the origin.   
\begin{prop}\label{compauto}
Assume that  $A\in \mathcal{M}_{2\times 2}(\mathbb{R})$  given in \eqref{Gensystem10}	is an excitable matrix.  For $j\in \mathbb{N}^{\ast},$ fix $\lambda_{j}>0$  given by \eqref{Gensyst7} and suppose $D=\text{diag}[d_{1},d_{2}],$ with $d_{i}>0,$ $i=1,2.$
Then,	the eigenvalues $\lambda_{1j}$ and $\lambda_{2j}$ of the matrix $B_{j}$ given in \eqref{ma2} have the following behavior:
	
	\begin{itemize}
		\item[i)]$\lambda_{1j}>0,$ and $\lambda_{2j}<0,$ if $(d_{1},d_{2})\in R_{u}.$
		
		\item[ii)] $\lambda_{1j}=0,$ and $\lambda_{2j}<0,$ if $(d_{1},d_{2})\in H_{\lambda_{j}}.$
		
		\item[iii)]$\mathcal{R}e(\lambda_{1j})<0,$ and $\mathcal{R}e(\lambda_{2j})<0,$ if $(d_{1},d_{2})\notin (R_{u}\cup H_{\lambda_{j}}).$ 
	\end{itemize}
\end{prop}
\begin{proof}
The roots of $P_{B_{j}}(\rho)$ are given by
	
	$$\lambda_{ij}=\frac{Tr(B_{j})\pm\sqrt{[Tr(B_{j})]^{2}-4\cdot \det(B_{j})}}{2},\;\; i=1,2.$$
Thus, if $(d_{1},d_{2})\in R_{u},$ then $\det(B_{j})<0$ and 
	$$\lambda_{1j}=\frac{Tr(B_{j})+\sqrt{[Tr(B_{j})]^{2}-4\cdot \det(B_{j})}}{2}>0\;\;\text{y}\;\;\lambda_{2j}=\frac{Tr(B_{j})-\sqrt{[Tr(B_{j})]^{2}-4\cdot \det(B_{j})}}{2}<0.$$

If $(d_{1},d_{2})\in H_{\lambda_{j}},$ then $\det(B_{j})=0.$ Therefore
	
	$$\lambda_{1j}=0\;\;\text{and}\;\;\lambda_{2j}=\text{Tr}(B_{j})<0.$$

	Finally, if $(d_{1},d_{2})\notin R_{u}\cup H_{\lambda_{j}},$ then $\det(B_{j})>0.$  The roots of $P_{B_{j}}(\rho)$ can be real or complex conjugate. If $[Tr(B_{j})]^{2}-4\cdot \det(B_{j})\geq 0,$ then the real roots satisfy 
	$$
	\lambda_{1j}=\frac{Tr(B_{j})+\sqrt{[Tr(B_{j})]^{2}-4\cdot \det(B_{j})}}{2}<0;\;\;\lambda_{2j}=\frac{Tr(B_{j})-\sqrt{[Tr(B_{j})]^{2}-4\cdot \det(B_{j})}}{2}<0.
$$ 
Note that, by continuity, $\lambda_{1j}$ and $\lambda_{2j}$ remain real roots of $P_{B_{j}}(\rho)$ provided that $(d_{1},d_{2})$ is near $H_{\lambda_{j}}.$
On the other hand, if  $[Tr(B_{j})]^{2}-4\cdot \det(B_{j})< 0,$ then the roots are complex conjugate and $ \displaystyle{\mathcal{R}e(\lambda_{1j})= \mathcal{R}e(\lambda_{2j})= \frac{\text{Tr}(B_{j})}{2}<0}.$
\end{proof}

\begin{rem}\label{eigenvectorsform}
In the second case of Proposition \ref{compauto}. i.e., when $(d_{1},d_{2})\in H_{\lambda_{j}}$ ($\det(B_{j})=0$), it is easy to show that
\begin{equation}\label{eigenvectorsform1}
v_{1j}= \left(
\begin{array}{c}
	1	\\
	-\frac{(a_{11}-\lambda_{j}d_{1})}{a_{12}}\\
\end{array}
\right) \;\text{and}\; v_{2l}= \left(
\begin{array}{c}
	1\\
	-\frac{(a_{11}-\lambda_{j}d_{1}-\lambda_{2j})}{a_{12}}\\
\end{array}
\right)
\end{equation}
or
\begin{equation}\label{eigenvectorsform2}
v_{1j}= \left(
\begin{array}{c}
	-\frac{(a_{22}-\lambda_{j}d_{2})}{a_{21}}	\\
	1\\
\end{array}
\right) \;\text{and}\; v_{2l}= \left(
\begin{array}{c}
	-\frac{(a_{22}-\lambda_{j}d_{2}-\lambda_{2j})}{a_{21}}	\\
1\\
\end{array}
\right)
\end{equation}
are the eigenvectors corresponding to eigenvalues
 $\lambda_{1j}=0,$ and $\lambda_{2j}=\text{Tr}(B_{j})<0,$ respectively. Since  $\det(B_{j})=0$ and $a_{12}, a_{21}$ have opposite signs then we infer  $a_{11}-\lambda_{j}d_{1} \neq 0$ and $a_{22}-\lambda_{j}d_{2} \neq 0.$ On the other hand, since $\lambda_{2j}=\text{Tr}(B_{j})$ one can show that $\det(B_{j}-\lambda_{2j}I_{2\times 2})=0,$ $a_{11}-\lambda_{j}d_{1} -\lambda_{2j}\neq 0$ and $a_{22}-\lambda_{j}d_{2} -\lambda_{2j}\neq 0.$ Other representations can be selected for the eigenvectors $v_{1j}$ and $v_{2j}$ but all representations will satisfy the following property:
 \begin{equation}\label{eigenpro}
 	v_{1j}=\left(
 	\begin{array}{c}
 		\xi_{1}	\\
 		\xi_{2}\\
 	\end{array}
 	\right)\;\text{with}\; \xi_{1} \neq 0,\; \xi_{2} \neq 0 \; \text{and}\; v_{2j} \; \text{is not parallel to}
 	\left(
 	\begin{array}{c}
 		\xi_{1}	\\
 		0\\
 	\end{array}
 	\right)\;\text{nor}\;
 	\left(
 	\begin{array}{c}
 		0	\\
 		\xi_{2}\\
 	\end{array}
 	\right).
 \end{equation}
Property \eqref{eigenpro} will be crucial to prove patterns existence via Turing bifurcation for system \eqref{Gensyst1} with $k=2$ (see Theorem \ref{Gensystem12} and \ref{Gensystem18}).
\end{rem}

	From Proposition \ref{compauto}, we observe that a real simple root of polynomial \eqref{charpolik=2} may only cross the imaginary axis at the origin in which case $\text{det}(B_{j})= 0.$ This together with the excitability of the Jacobian matrix $A$ (see \eqref{Gensystem10}) are necessary conditions for the change of stability of the matrix $B_{j}$ (see \eqref{ma2}), and therefore for the associated bifurcation to occur, we describe this in the following Section.

 \section{Patterns formation for the reaction-diffusion system \eqref{Gensyst1} with $k=2$ via Turing bifurcation}\label{Turing Bufurcation}
 In this section we will show that the diffusion-driven instability phenomenon gives rise to non-homogeneous steady-state solutions of \eqref{Gensyst1} with $k=2$ that bifurcate from the uniform stationary solution $w_{0}=(w_{0}^{1}, w_{0}^{2}) \in \mathring{\Lambda}_{2}.$ Specifically, we establish sufficient conditions for a Turing bifurcation to occur in the reaction-diffusion system \eqref{Gensyst1} with $k=2$ by using  \cite[Theorem 13.5]{Smoller}, which is a bifurcation from a  simple eigenvalue theorem for operators in Banach spaces. For this, we define
 \begin{equation}\label{AltX}
\tilde{X}:=\{w\in C^{2}(\overline{\Omega}, \mathbb{R}) \times C^{2}(\overline{\Omega}, \mathbb{R}): \; \displaystyle{\frac{\partial w}{\partial \eta}}(x,t)=0, \;t>0, \; x\in \partial \Omega \}
\end{equation}
with the usual supremum norm involving the first and second derivatives
$$\|w\|_{\tilde{X}}:=\sum_{i=1}^{2} \frac{1}{i!}\|w^{(i)}(x)\|_{\infty},$$
and $Y:=C(\overline{\Omega}, \mathbb{R}) \times C(\overline{\Omega}, \mathbb{R}) $ with the usual supremum norm. However, when choosing the subspace $Z$ in \cite[Theorem 13.5]{Smoller},  we will use the orthogonality induced by the scalar product
$$\left\langle v,w\right\rangle=\int\limits_{\overline{\Omega}}(v_{1}(x) w_{1}(x)+v_{2}(x) w_{2}(x))dx,$$
where $v=(v_{1},v_{2})$ and $w=(w_{1},w_{2}).$ Now we are able to prove our first criteria which consider $d_{1}$ as bifurcation parameter. 

\begin{proof}[Proof of Theorem \ref{Gensystem12}]
	Setting $u=w-w_{0},$ where $w_{0}$ is a non-trivial homogeneous steady state solution of  \eqref{Gensyst1} with $k=2,$ we get
$$	\left \{
\begin{array}{llll}
	\displaystyle{\frac{\partial u(x,t)}{\partial t}} = D \displaystyle{\Delta u(x,t)}+Au(x,t)+ Q(u(x,t)),\;\;\;x\in\Omega,\;\;t>0,\\
	\displaystyle{\frac{\partial u}{\partial \eta}}(x,t)=0,\;\;\; x \in \partial \Omega,\;\;t>0,	
\end{array}
\right.	$$
where $A$ is the Jacobian matrix at $\omega_{0}\in \mathring{\Lambda}_{2}$ given in \eqref{Gensystem10} satisfying the hypothesis of the theorem and $Q(u)=G(w_{0}+u)-Au.$\\

For any nonhomogeneous stationary solution $w$ of \eqref{Gensyst1} with $k=2,$ $u=w-w_{0}$ satisfies the elliptic equation 
\begin{equation}\label{ellip1}
	\left \{
\begin{array}{llll}
 D \displaystyle{\Delta u(x,t)}+Au(x,t)+ Q(u(x,t))=0,\;\;\;x\in\Omega,\;\;t>0,\\
	\displaystyle{\frac{\partial u}{\partial \eta}}(x,t)=0,\;\;\; x \in \partial \Omega,\;\;t>0.	
\end{array}
\right.	
\end{equation}
Taking into account this observation, we define the function $\tilde{F}$ and the linear operator $L_{0}$  as follows: 
\begin{equation}\label{deLo}
	\begin{split}
		\tilde{F}:\mathbb{R}^{+}\times \tilde{X}\to Y \qquad \qquad&\qquad \qquad \qquad \quad L_{0}:\tilde{X} \to Y\\
		\tilde{F}(d_{1},u)= D \displaystyle{\Delta u}+Au+ Q(u)\qquad &\qquad L_{0}=D_{2}\tilde{F}(d_{1}^{\ast},0)= \left(
		\begin{array}{cc}
			d_{1}^{\ast}& 0\\
			0 &d_{2}\\
		\end{array}
		\right)\Delta +A,
	\end{split}
\end{equation}
where $\tilde{X},\; Y$ are defined in \eqref{AltX}, $d_{1}$ is the diffusion coefficient of the susceptible class, $d_{2}$ is a positive real number satisfying condition \eqref{Gensystem14},  and the differential with respect to the $u$ variable at $(d_{1}^{\ast},0)$ is taken in the Fréchet sense (see \cite[Definition 13.1]{Smoller}). From \eqref{Gensyst7} we infer that, if $v_{ij},$ with $i=1,2$ and $j=0,1,2...,$ denote the eigenvectors of $B_{j}=A-\lambda_{j}D$ evaluated at $d_{1}=d_{1}^{\ast}$ with respective eigenvalues $\lambda_{ij},$ then  the spectrum of the linear operator $L_{0}$  is given by the eigenvalues $\lambda_{ij}$  with $\psi_{ij}:=v_{ij}\phi_{j}$ as respective eigenfunctions. In fact,
$$\lambda_{ij}v_{ij}=B_{j}v_{ij} \Longleftrightarrow \lambda_{ij}v_{ij}\phi_{j}= B_{j}v_{ij}\phi_{j}=L_{0}v_{ij}\phi_{j},$$
where $\phi_{j}$ denote the $j$th eigenfunction of $-\Delta$ on $\Omega$ with no flux-boundary conditions (see \eqref{Gensyst7}).\\

A simple algebra calculation shows that, in the positive plane $d_{1}d_{2},$ each Hyperbola $H_{\lambda_{j}}$  defined by \eqref{Gensystem11} intersects the others exactly once and since  $d_{2}'<d_{2}<d_{2}'',$ then  $l\geq 2$ is the unique natural number such that $(d_{1}^{\ast}, d_{2})$ belongs to the hyperbola $H_{\lambda_{l}}$, see Figure \ref{Bifurca22}. Furthermore, from  commentaries in subsection \ref{eigenbeha} and Proposition \ref{compauto} we have that $\text{det}(B_{j})>0$ for $j \neq l$ and $\text{det}(B_{j})=0$ just for $j=l.$ Hence, for $i=1,2$ and $j=0,1,2,...,l-1,l+1,...$ all the eigenvalues $\lambda_{ij}$ have negative real part. For $j=l,$ the matrix $B_{l}$ given in  \eqref{Gensystem16} have one eigenvalue, say $\lambda_{1l},$ equal to zero and the other one is negative, i.e. $\lambda_{2l}<0.$ Furthermore, the corresponding eigenvectors denoted by $v_{1l}=\left(
\begin{array}{c}
	\xi_{1}	\\
	\xi_{2}\\
\end{array}
\right)$ and $v_{2l}$ satisfy the property \eqref{eigenpro}. Thus, the eigenfunction of the linear operator $L_{0}$ corresponding to $\lambda_{1l}=0$ is given by $\psi_{1l}=v_{1l}\phi_{l}$ which is a non-uniform stationary solution of the linearized system \eqref{Gensyst1} with $k=2,$ that is,
\begin{equation}\label{Ellipt}
	\left \{
	\begin{array}{llll}
\left(
\begin{array}{cc}
	d_{1}^{\ast}& 0\\
	0 &d_{2}\\
\end{array}
\right) \displaystyle{\Delta \psi_{1l}(x)}+ A\psi_{1l}(x)=0,\;\;\;x\in\Omega,\\
		\displaystyle{\frac{\partial \psi_{1l}}{\partial \eta}}(x)=0,\;\;\; x \in \partial \Omega.	
	\end{array}
	\right.
\end{equation}

\begin{figure}[h!]
	\includegraphics[width=12cm]{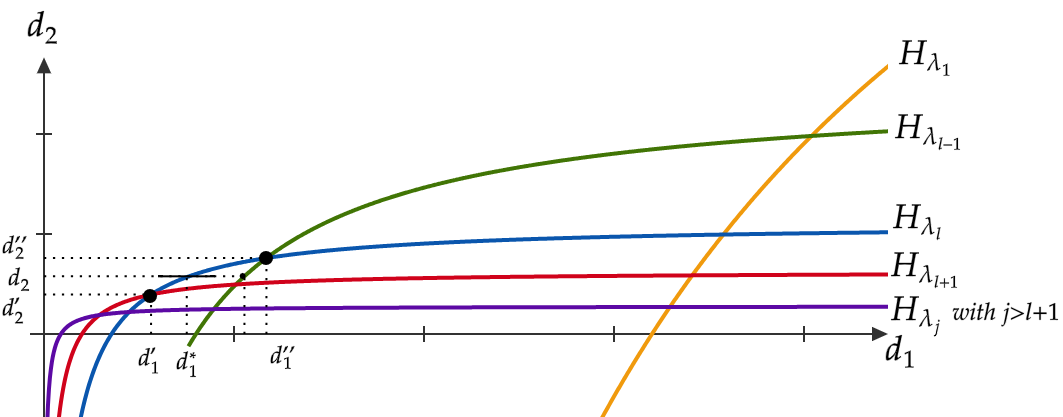}
	\caption{Turing Bifurcation: when $d_{2}'<d_{2}<d_{2}'',$ the uniform steady state solution $w_{0}\in \mathring{\Lambda}_{2}$ of \eqref{Gensyst1} with $k=2$ undergoes a Turing bifurcation at $d_{1}=d_{1}^{\ast}.$  Here, $\lambda_{j}$ denote the $j$th eigenvalue of $-\Delta$ on $\Omega$ with no flux-boundary conditions (see \eqref{Gensyst7}) and the Hyperbola $H_{\lambda_{j}},$ $j=1,2,...,l-1,l, l+1,...$ is defined in \eqref{Gensystem11} (see Figure \eqref{casoa22} part $b).$)}\label{Bifurca22}
\end{figure}
Therefore, the null-subspace $\mathcal{N}(L_{0})$ of operator $L_{0}$ is one-dimensional spanned by $\psi_{1l}.$ The Range $\mathcal{R}(L_{0})$ of this operator is given by the relation 
$${\normalsize
	\begin{split}
		\mathcal{R}(L_{0}):=\{
		z\in [C(\Omega,\mathbb{R})]^{2}:\;&z'\text{s}\; \text{Fourier expansion does not contain the term}\; \phi_{l} \} \cup \{v_{2l}\phi_{l}\}
\end{split}}
$$
because the orthogonality and completeness of the system $\phi_{j},$ $j=0,1,2,...$ obtained by solving the eigenvalue problem \eqref{Gensyst7}.
So, the codimension of $\mathcal{R}(L_{0})$ is one and conditions $i)$ and $ii)$ of \cite[Theorem 13.5]{Smoller} are satisfied. It still remain to verify condition $iii).$ Let $L_{1}$ be  the differential Fréchet operator of $D_{2}\tilde{F}$ with respect to the $d_{1}$ variable at $(d_{1}^{\ast},0).$ Then
$${\normalsize 
	\begin{split}
 L_{1}:\tilde{X}& \to Y\\
 L_{1}=D_{1}D_{2}\tilde{F}(d_{1}^{\ast},0)&
 = \left(
		\begin{array}{cc}
			1& 0\\
			0 &0\\
		\end{array}
		\right)\Delta,
	\end{split}
}$$
\begin{equation}\label{defLl}
	L_{1}\psi_{1l}=\left(
	\begin{array}{cc}
		1& 0\\
		0 &0\\
	\end{array}
	\right)\Delta v_{1l}\phi_{l}
	= \left(
	\begin{array}{cc}
		1& 0\\
		0 &0\\
	\end{array}
	\right) 
	\left(
	\begin{array}{c}
		\xi_{1}	\\
		\xi_{2}\\
	\end{array}
	\right) \Delta \phi_{l}
	=-\lambda_{l}
	 \left(
	\begin{array}{c}
		\xi_{1}	\\
		0\\
	\end{array}
	\right) \phi_{l},
\end{equation}
with  $\displaystyle{ \left(
	\begin{array}{c}
		\xi_{1}	\\
		0\\
	\end{array}
	\right) }$ not being parallel to $v_{2l}$ (see \eqref{eigenpro}), and
$$\left\langle \psi_{1l},L_{1}\psi_{1l} \right\rangle=\int\limits_{\overline{\Omega}}-\lambda_{l} \xi_{1}^{2} \phi_{l}^{2}(x)dx \neq 0$$
because $\xi_{1}\neq 0.$ Thus, $L_{1}\psi_{1l} \notin \mathcal{R}(L_{0})$ and condition $iii)$ of \cite[Theorem 13.5]{Smoller} is satisfied. Hence, by choosing $Z=\mathcal{R}(L_{0}),$ we conclude that there exists $\zeta>0$ and a $C^{1}-$curve $(d_{1},z):(-\zeta, \zeta)\to \mathbb{R}^{+} \times Z$ with $d_{1}(0)=d_{1}^{\ast},$ $z(0)=0$ and $\tilde{F}(d_{1}(s),s\cdot[\psi_{1l}+z(s)])=0.$ Therefore,
$$
u(x,s)=s\cdot v_{1l}\phi_{l}(x)+s\cdot z(x,s)
$$
is a solution of the elliptic equation \eqref{ellip1} where $d_{1}\equiv d_{1}(s),$  $s\in (-\zeta, \zeta)$ and $v_{1l}$ is given by \ref{eigenvectorsform1} with $j=l$ and $d_{1}=d_{1}^{\ast}.$ Finally, taking into account that $u=w-w_{0}$ we obtain that $w(x,s)$ defined as in \eqref{Gensystem15} are non-uniform stationary solutions of \eqref{Gensyst1} with $k=2.$   Since $s$ is considered to be small, we call this solution a \textit{small amplitude pattern}. Therefore, at $d_{1}=d_{1}^{\ast},$ the uniform steady state solution $w_{0}$ undergoes a Turing bifurcation. This proves the Theorem.
\end{proof}

Now we show our second criterion.
The proof is similar to that of Theorem \ref{Gensystem12} (see also \cite[Theorem 3.2]{Cavani and Farkas II} and \cite[Theorem 2]{Lizana and Marin}). So, we bring only the necessary changes. In this case we determine nonhomogeneous stationary solutions of \eqref{Gensyst1} with $k=2$ considering $d_{2}$ as a bifurcation parameter.

 \begin{proof}[Proof of Theorem \ref{Gensystem18}]
Here, we define the function $\tilde{F}$ and the linear operator $L_{0}$  as follows: 
  $${\normalsize 
  	\begin{split}
  		\tilde{F}:\mathbb{R}^{+}\times \tilde{X}\to Y \qquad \qquad&\qquad \qquad \qquad \quad L_{0}:\tilde{X} \to Y\\
  		\tilde{F}(d_{2},u)= D \displaystyle{\Delta u}+Au+ Q(u)\qquad &\qquad L_{0}=D_{2}\tilde{F}(d_{2}^{\ast},0)= \left(
  		\begin{array}{cc}
  			d_{1}& 0\\
  			0 &d_{2}^{\ast}\\
  		\end{array}
  		\right)\Delta +A,
  	\end{split}
  }$$
  where $d_{2}$ is the diffusion coefficient of the susceptible class and $d_{1}$ is a positive real number satisfying condition \eqref{Gensystem19}. Therefore,
  $${\normalsize 
  	\begin{split}
  		L_{1}:\tilde{X}& \to Y\\
  		L_{1}=D_{1}D_{2}\tilde{F}(d_{2}^{\ast},0)&
  		= \left(
  		\begin{array}{cc}
  			0& 0\\
  			0 &1\\
  		\end{array}
  		\right)\Delta,
  	\end{split}
  }$$
and the results established in Theorem \ref{Gensystem18} follows.
 \end{proof}

\section{Stability of Bifurcating Solutions}\label{stab} 
 In this section we analyze the stability of the one-parameter family  of non-uniform stationary solution  of system \eqref{Gensyst1} with $k=2$ that arise from the bifurcation of the homogeneous steady state $w_{0}.$ Our first result study the stability of the bifurcating nonhomogeneous stationary solution \eqref{Gensystem15}.

 \begin{proof}[Proof of Proposition \ref{stability}]
 	In the proof of Theorem \ref{Gensystem12}, we have shown that $\lambda_{1l}=0$ is a $L_{1}-$simple eigenvalue of $L_{0}$ (see \eqref{deLo}-\eqref{defLl}) with eigenfunction $\psi_{1l}$ (see \cite[Definition 13.6]{Smoller}). In particular, $L_{0}\psi_{1l}=\lambda_{1l} \psi_{1l}.$ 
 	
 	On the other hand, note that if $d_{1}\in (d_{1}^{\ast}-\epsilon, d_{1}^{\ast}+\epsilon),$ $s\in(-\zeta, \zeta),$ and $|\epsilon|,$ $|\zeta|$ are small enough, then the operator $D_{2}\tilde{F}(d_{1},0)$ and
 	$D_{2}\tilde{F}(d_{1}(s),s\psi_{1l} + s z(x,s))$
 	 are close to $L_{0}.$ Applying \cite[Lemma 13.7]{Smoller} for both operators, we obtain that there exists smooth functions 
$${\normalsize
	\begin{split}
	(d_{1}^{\ast}-\epsilon, d_{1}^{\ast}+\epsilon)	\to \mathbb{R} \times \tilde{X} \qquad \qquad \qquad \qquad&  \qquad  (-\zeta, \zeta) \to \mathbb{R} \times \tilde{X}\\
		d_{1} \longmapsto (\rho(d_{1}), \psi_{c}(d_{1}))\qquad \qquad &\qquad 	\qquad	s \longmapsto (\gamma(s), \psi_{b}(s)),
	\end{split}}
$$	 
such that 
$${\normalsize
	\begin{split}
D_{2}\tilde{F}(d_{1},0)\psi_{c}(d_{1})& = \rho(d_{1}) \psi_{c}(d_{1}),\\
D_{2}\tilde{F}(d_{1}(s),s\psi_{1l}(s)+sz(x,s))\psi_{b}(d_{1})& = \gamma(s) \psi_{b}(s),\\
	\end{split}
}$$
 	and
$(\rho(d_{1}^{\ast}),\psi_{c}(d_{1}^{\ast}))=(0,\psi_{1l})=(\gamma(0), \psi_{b}(0)).$ 	Note that, when applying  \cite[Lemma 13.7]{Smoller} for  operator $D_{2}\tilde{F}(d_{1},0)$ we have setting $\rho(d_{1})=\eta(D_{2}\tilde{F}(d_{1},0))$ and 
$\psi_{c}(d_{1})=\text{w}(D_{2}\tilde{F}(d_{1},0)).$ Similarly, when applying  \cite[Lemma 13.7]{Smoller} for  operator $D_{2}\tilde{F}(d_{1}(s),s\psi_{1l}(s)+sz(x,s))$ we have setting $\gamma(s)=\eta(D_{2}\tilde{F}(d_{1}(s),s\psi_{1l}(s)+sz(x,s)))$ and $\psi_{b}(s)=\text{w}(D_{2}\tilde{F}(d_{1}(s),s\psi_{1l}(s)+sz(x,s))).$\\

It is easy to show that 
\begin{equation}\label{derdifzero}
\rho'(d_{1}^{\ast})\neq 0.
\end{equation}
In fact, we know that $\rho(d_{1})$ satisfies the equation 
$$ 	\rho^{2}(d_{1})-\text{Tr}(B_{l}) \rho(d_{1}) +\text{det}(B_{l})=0,$$
(see \eqref{ma2}-\eqref{charpolik=2}). Applying implicit differentiation  with respect to $d_{1}$ on this equation,  we have
$$\rho'(d_{1})=\frac{\lambda_{l}a_{22}-\lambda_{l}^{2} d_{2}-\lambda_{l} \rho(d_{1})}{2 \rho(d_{1})-\text{Tr}(B_{l})}.$$ 
Therefore, 
$$\rho'(d_{1}^{\ast})=\frac{\lambda_{l}(\lambda_{l} d_{2}-a_{22})}{\text{Tr}(A)-\lambda_{l}(d_{1}^{\ast}+d_{2})}.$$
From the fact that $a_{22}>0$ and $0<d_{2}<\frac{a_{22}}{\lambda_{l}}$ (see part $b)$ of Figure \ref{casoa22}), we infer that $\rho'(d_{1}^{\ast})> 0$ and \eqref{derdifzero} holds.

From hypothesis, $\gamma(s)\neq 0$ for $s$ close to $0,$  then a direct application of the Crandall-Rabinowitz Theorem \cite[Theorem 1.16]{Crandal Rabinowitz} (see also \cite[Theorem 13.8, pp. 180]{Smoller}) implies that
\begin{equation}\label{limite}
\lim\limits_{s\to 0} \frac{s d_{1}'(s) \rho'(d_{1}^{\ast})}{\gamma(s)}=-1.
\end{equation}

Finally, we use \eqref{limite} to determine the sign of $\gamma(s)$ and  consequently the stability of the bifurcating solution \eqref{Gensystem15}. In fact, from hypothesis $d_{1}'(0)\neq 0$ then, without loss of generality, one can assume that $d_{1}'(0)>0.$ Hence, by continuity, we infer that $d_{1}'(s)>0$ for $|s|$ small enough. From \eqref{limite} it follows that $\gamma(s)<0$ for $s>0$  small enough ($d_{1}(s)>d_{1}^{\ast}$ because $d_{1}(s)$ is an increasing function) and therefore the bifurcating solution \eqref{Gensystem15} is asymptotically stable. On the other hand, if $s<0$ for small enough $s$ ($d_{1}(s)<d_{1}^{\ast}$) then \eqref{limite} implies that $\gamma(s)>0$ and the bifurcating nonhomogeneous stationary solution \eqref{Gensystem15} is unstable. The case $d_{1}'(0)<0$ can be analyzed similarly.
 \end{proof}
 
 Using similar arguments as in the  proof of Proposition \ref{stability} we can show Proposition  \ref{stability2}.
 
 \section{Applications}\label{app} 
 
 The interaction of at least two species with considerably different diffusion coefficients can give rise to pattern formation (see \cite{Malshow Petrovskii and Venturino} and the references therein). The particular example of Segel and Jackson in \cite{Segel and Jackson}
 applied Turing's idea to a problem in population dynamics: the Turing instability in the prey-predator interaction of algae and herbivorous copepods with higher herbivore motility. Specifically, as exposed by the authors in \cite{Segel and Jackson}, in situations where uneven geographic distribution of predator and prey would be mutually advantageous. For instance, if the prey were capable of some sort of cooperation so that the number of offspring per prey individual was an increasing function of prey density at first (of course, one would expect that the birth rate would be a decreasing function of prey density at higher density levels). It would appear to be mutually beneficial if the predators concentrated in certain areas, letting the prey population rise outside the areas of predator concentration. At higher population levels, the prey’s ability to cooperate would allow them to reproduce faster. The predators would partially benefit from this, since some of the larger prey population would “diffuse” into the concentrations of predator. \\

 Inspired by the works in \cite{Segel and Jackson, Lizana and Marin, Cavani and Farkas II}, in this paper we intend to apply Criterion I and II (see Theorems \ref{Gensystem12} and \ref{Gensystem18}) to  study the existence of patterns for the following reaction-diffusion  predator-prey model with variable mortality and Hollyn's type II functional response.
 \begin{equation}\label{FunctionResp}
 	\left \{
 	\begin{array}{llll}
 		\displaystyle{\frac{\partial u}{\partial t}} = d_{1} \displaystyle{\Delta u}+ Au\left(1-\frac{u}{K}\right)-a \frac{uv}{1+Eu},\;\;\;x\in\Omega,\;\;t>0,\\
 		\\
 		\displaystyle{\frac{\partial v}{\partial t}} = d_{2} \displaystyle{\Delta v}-\left[\delta + \frac{\gamma - \delta}{1+v}\right]v+b\frac{uv}{1+Eu},\;\;\;x\in\Omega,\;\;t>0,\\
 	\end{array}
 	\right.
 \end{equation}
 subject to homogeneous Newmann boundary conditions
 \begin{equation}\label{NewmanCond}
 	\displaystyle{\frac{\partial u}{\partial \eta}}(t,x)= \displaystyle{\frac{\partial v}{\partial \eta}}(t,x)=0,\;\;\; x \in \partial \Omega,\;\;t>0,
 \end{equation}
 and initial data
 \begin{equation}\label{InitialCond}
 	u(0,x)=\varphi_{1}(x)\geq0;\;\;v(0,x)=\varphi_{2}(x)\geq0,\;\;x\in\Omega,
 \end{equation}
 where $A,$ $K,$ $a,$ $E,$ $\delta,$ $\gamma,$ $b$ are positive constants and $u(x,t),$ $v(x,t)$ represent the population density of prey and predator at $x \in \Omega$ and time $t,$ respectively. $d_{i}>0,$ $i=1,2$ are the diffusive coefficients of the prey and predator respectively and 
 $\Omega \subset \mathbb{R}^{n}$ as usual is a bounded, open and connected set.
 The prey grows with intrinsic growth rate $A$ and carrying capacity $K$ in the absence of predation. The parameter $a$ represent the conversion rate with respect to the prey. The predator consumes the prey with functional response $\frac{uv}{1+Eu}$ and satiation coefficient or conversion rate $b;$ the specific mortality of predators in absence of prey 
 \begin{equation}\label{mort}
 	M(v)=\delta + \frac{\gamma - \delta}{1+v}= \frac{\gamma + \delta v}{1+v},\qquad 0< \gamma < \delta,
 \end{equation}
 depends on the quantity of predators; $\gamma$ is the mortality at low density and $\delta$ is the maximal mortality; the natural assumption is $\gamma< \delta.$ The advantage of this model over the most often used models is that here the predator mortality is neither a constant nor an unbounded function, still it is increasing with quantity. 
 \begin{figure}[h!]
 	\includegraphics[width=6cm]{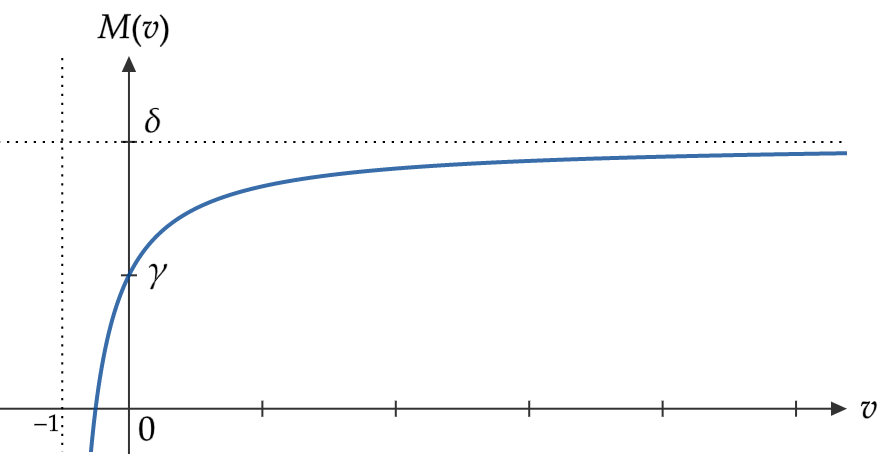}
 	\caption{Graph of the predator mortality function $M(v).$}\label{PMortal}
 \end{figure}
 The variable mortality \eqref{mort} was introduced in \cite{Cavani and Farkas I, Cavani and Farkas II, Lida and Lizana} where the authors studied bifurcations (see \cite{Cavani and Farkas I, Cavani and Farkas II}) and found homoclinic orbits in a predator-prey model. In their work, Cavani and Farkas \cite{Cavani and Farkas II}  considered a functional response of Holling type and $\Omega$ as an closed interval of $\mathbb{R}.$ It is known that the patterns produced by Turing's mechanism can be sensitive to domain's shape \cite{Bunow et al}, so in this section we investigate this issue for system \eqref{FunctionResp}-\eqref{NewmanCond} by performing numerical simulations.
 Setting  $z=(u,v),$ $\varphi=(\varphi_{1},\varphi_{2}),$ $D=\text{diag}[d_{1},d_{2}]$ and $P=(P_{1},P_{2}),$ where
$${\normalsize 
 	\begin{split}
 		P_{1}(z)=P_{1}(u,v):=& Au\left(1-\frac{u}{K}\right)-a \frac{uv}{1+Eu},\\
 		P_{2}(z)=P_{2}(u,v):=&-\left[\delta + \frac{\gamma - \delta}{1+v}\right]v+b\frac{uv}{1+Eu},
 	\end{split}}
$$
 the system \eqref{FunctionResp}-\eqref{NewmanCond}-\eqref{InitialCond} can be rewritten as
 
 \begin{equation}\label{Gensystem}
 	\left \{
 	\begin{array}{llll}
 		\displaystyle{\frac{\partial z(x,t)}{\partial t}} = D \displaystyle{\Delta z(x,t)}+ P(z),\;\;\;x\in\Omega,\;\;t>0,\\
 		\displaystyle{\frac{\partial z}{\partial \eta}}(x,t)=0,\;\;\; x \in \partial \Omega,\;\;t>0,		\\
 		z(x,0)=\varphi(x),\;\;x\in\Omega.
 	\end{array}
 	\right.
 \end{equation}

 \subsection{Well-posedness of system \eqref{Gensystem}}\label{wellpose} 
 
 Here we focus on the well-posedness of system \eqref{Gensystem}. In fact, a direct application of Theorem \ref{GenerWP} shows that the system \eqref{Gensystem} is biologically well-posed and its relevant dynamic is concentrated in $X_{\Lambda_{2}}.$ This is stated in the following result.
 \begin{cor}\label{existandunic}
 For each $\varphi=(\varphi_{1}, \varphi_{2}) \in X_{\Lambda_{2}}$ (see \eqref{Gensyst2}), the system \eqref{Gensystem} has a unique mild solution $z(t)=z(t, \varphi)\in X_{\Lambda_{2}}$ defined on $[0,\sigma)$ and a classical solution $z(x,t)\equiv [z(t)](x).$ Moreover, the set $X_{\Lambda_{2}}$ is positively invariant under the flow $\varPsi_{t}(\varphi)=z(t,\varphi)$ induced by \eqref{Gensystem}. 
 \end{cor}
 \begin{proof}
 Note that the nonlinear term $P$ present in \eqref{Gensystem} satisfies $P\in C^{2}( \mathring{\Lambda}_{2},  \mathbb{R}^{2}),$ 
 	and $P_{i}(z)\geq 0$ whenever $z=(z_{1},z_{2})=(u,v) \in  \Lambda_{2}$ and $z_{i}=0.$ Theorem \ref{GenerWP} gives the desired result.
 \end{proof}
 
 Next, we will show the parameter $\sigma$ obtained in Corollary \ref{existandunic} is equal to infinity. Therefore, all solutions of system \eqref{Gensystem} are defined for all $t\geq 0.$ This is established in the following lemma.

 \begin{lem}\label{compact}
 	Assume  $\varphi=(\varphi_{1}, \varphi_{2}) \in X_{\Lambda_{2}}.$
 	Let $z=(u,v)$ be any solution of \eqref{Gensystem}.  Then 
 	\begin{equation}\label{stm1}
 		\limsup_{t\to \infty} \max_{x\in \Omega} u(x,t)\leq K,
 	\end{equation}		
 	and if   $b-E\delta<0,$ then
 		\begin{equation}\label{stm2}
 			\limsup_{t\to \infty} \max_{x\in \Omega} v(x,t)\leq - \frac{E(\delta-\gamma)}{(b-\delta E)}.
 		\end{equation}
 \end{lem}
 
 \begin{proof}
 	From the first equation of system \eqref{Gensystem}, we infer that
 	$$ \displaystyle{\frac{\partial u}{\partial t}} \leq  d_{1} \displaystyle{\Delta u}+ Au\left(1-\frac{u}{K}\right),  $$
 	as long as $N$ is defined as a function of $t.$ This, give us the adequate function to apply the Comparison Principle  \cite[Chapter 7, $\S 3$]{Hal Smith} (see also  \cite[Chapter 11, $\S 11.1$ c, $\S 11.3$ a]{McOwen}). Let  $\rho$ be the solution of the ODE
 $$
 		\left \{
 		\begin{array}{llll}
 			\rho'(t)=A\rho(t)\left(1
 			-\frac{\rho(t)}{K}\right),\;\;t>0,\\
 			\rho(0)=\displaystyle{\max_{x \in \Omega} u(x,0).}
 		\end{array}
 		\right.
 $$
 	Note that $\rho(t)$ is well defined for $0<t<\infty$  and $\rho(t) \to K$ as $t\to \infty.$  From  Theorem 3.4 in \cite[Chapter 7, $\S 3$]{Hal Smith}   we have that $u(x,t)\leq \rho(t).$ Therefore, for any $\epsilon>0,$ there exists $T_{\epsilon}>0$ such that
 	\begin{equation}\label{cotasup}
 		u(x,t)\leq \rho(t)<K+\epsilon,\;\;(x,t)\in \overline{\Omega} \times [T_{\epsilon},+\infty).
 	\end{equation}
 	On the other hand,  if  $b-E\delta<0.$ Then, from \eqref{mort} and the second equation of system \eqref{Gensystem} we infer
 	$${\normalsize
 		\begin{split}
 			\displaystyle{\frac{\partial v}{\partial t}} &= d_{2} \displaystyle{\Delta v}-\left[\delta + \frac{\gamma - \delta}{1+v}\right]v+\frac{b}{E}\frac{(Eu)v}{(1+Eu)}
 			\leq  d_{2} \displaystyle{\Delta v}+\frac{(b-\delta E)v}{E} +(\delta-\gamma),
 	\end{split}}
 	$$
 for any $ x\in \Omega$ and $t\geq 0.$ Let  $s$ be the solution of the ODE
 $$
 		\left \{
 		\begin{array}{llll}
 			s'(t)=\displaystyle{\frac{(b-\delta E)}{E}s(t)} +(\delta-\gamma),\;\;t>0,\\
 			s(0)=\displaystyle{\max_{x \in \Omega} v(x,0).}
 		\end{array}
 		\right.
$$
 	Therefore, 
 	$$s(t)=-\frac{E(\delta-\gamma)}{b-\delta E}+\left(s(0)+\frac{E(\delta-\gamma)}{b-\delta E}\right)e^{\frac{(b-\delta E)t}{E}}\leq 
 	-\frac{E(\delta-\gamma)}{b-\delta E}+ s(0)e^{\frac{(b-\delta E)t}{E}}.$$
 	Again, the Comparison Principle implies
 	\begin{equation}\label{csupv}
 		v(x,t)\leq s(t) \leq 
 		-\frac{E(\delta-\gamma)}{b-\delta E}+ s(0)e^{\frac{(b-\delta E)t}{E}}, \;\;(x,t)\in \overline{\Omega} \times [0,+\infty),
 	\end{equation}
 	as long as $v$ is defined as a function of $t.$
 	
 	\vskip.2cm

 	From \eqref{cotasup}, \eqref{csupv} and part $iii)$ of Theorem \ref{GenerWP},  we conclude  that the solution $z(x,t)=(u(x,t),v(x,t))$ of the system \eqref{Gensystem} is defined for all $t\geq 0$ and satisfies the  estimates \eqref{stm1}--\eqref{stm2}. This proves the Lemma.
 \end{proof}
 
 \begin{rem}
 	Lemma \ref{compact} implies that the relevant dynamic of the system \eqref{Gensystem} is concentrated in a compact set of the space $X_{\Lambda_{2}}.$ 
 \end{rem}
 
   \subsection{Location and asymptotic stability analysis of equilibria for predator-prey model \eqref{FunctionResp}-\eqref{NewmanCond}.}\label{Applications}
   In this subsection we will study system \eqref{FunctionResp} without diffusion. Particularly, we will focus our attention on the existence of equilibria and their local stability. This information will be crucial in the next section where we study the effect of the diffusion parameters on the stability of the steady state. For this we consider the kinetic system associated to  
\eqref{FunctionResp},
  \begin{equation}\label{KineHollinhII}
 	\left \{
 	\begin{array}{llll}
 		\displaystyle{u' =  Au\left(1-\frac{u}{K}\right)-a \frac{uv}{1+Eu}},\\
 		\\
 	\displaystyle{	v' = -\left[\delta + \frac{\gamma - \delta}{1+v}\right]v+b\frac{uv}{1+Eu}}.\\
 	\end{array}
 	\right.
 \end{equation}
where $A,$ $K,$ $a,$ $E,$ $\delta,$ $\gamma,$ $b$ are positive constants.
Note that the positive quadrant of the phase plane $\Lambda_{2}:=\{(u,v)\in \mathbb{R}^{2}:u\geq 0, v\geq 0\}$ is positively invariant with respect to the system \eqref{KineHollinhII}.
It is easy to show that $(0,0)$ and $(K,0)$ are equilibrium points.  A simple linear stability analysis shows that $(0,0)$ is a saddle point and therefore it is always unstable. On the other hand, $(K,0)$ is local asymptotically stable if 
\begin{equation}\label{lasympsta}
	(b-\gamma E)K<\gamma,
\end{equation}
and unstable (a saddle point) if
\begin{equation}\label{unstsadd}
	(b-\gamma E)K > \gamma.
\end{equation}
Note that \eqref{unstsadd} implies  $b>E\gamma.$ Nevertheless, for reasonable parameter configurations we may establish the global stability of $(K,0).$

\begin{lem}\label{asympt}
	If  $\displaystyle{\delta > \gamma \geq \frac{b}{E}}$  then $(K,0)$ is global asymptotically stable with respect to the positive quadrant of the $uv-$plane.
\end{lem}
\begin{proof}
	Under assumptions of the lemma we have
	$${\normalsize
	\begin{split}
		v'&=-v\left[\frac{\gamma+\delta v}{1+v}-b \frac{u}{1+Eu}\right] 
		\leq -v \left[\gamma- \frac{b}{E} \frac{u}{(\frac{1}{E}+u)}\right]
		\leq -\frac{b}{E}v \left[1-  \frac{u}{(\frac{1}{E}+u)}\right],		
		\leq -Cv,		
\end{split}}$$
for some $C>0,$ since $u(t)$ is bounded for $t\in [0,+\infty).$ Hence, any solution of $v(t)$ corresponding to non-negative initial conditions tend to zero as $t \to \infty.$ Therefore, the omega limit set $\omega$ of every solution with positive initial conditions is contained in $\{(u,0)\in \mathbb{R}^{2}: u\geq 0\}.$ On the other hand, note that for $u>K,$ we have $u'<0.$ Thus, $\omega \subset \{(u,0)\in \mathbb{R}^{2}: 0\leq u \leq K\}.$ Since $(0,0)\notin \omega,$ $\omega \neq \emptyset,$ closed, and an invariant set we obtain that $\omega = (K,0).$ 
\end{proof}

\begin{rem}
	In particular, 	if $E>1,$ and $\delta> \gamma \geq b$  then $(K,0)$ is global asymptotically stable with respect to the positive quadrant of the $uv-$plane.
	Condition $E\gamma > b$ in Lemma \eqref{asympt} means that the minimal mortality of predator scaled by $E$ is high compared with the conversion rate.
\end{rem}

\begin{lem}\label{globsta}
If $\displaystyle{\delta \geq \frac{b}{E}>\gamma,} $ $a\geq 1$ and 
\begin{equation}\label{asymp}
(b-\gamma E)K\leq \gamma
\end{equation}
 then $(K,0)$ is global asymptotically stable	with respect to the positive quadrant of the $uv-$plane.
\end{lem}
\begin{proof}
	Assume first that inequality \eqref{asymp} is strict. Then, \eqref{lasympsta} holds and $(K,0)$ is local asymptotically stable. Using $b>\gamma E$ we infer that there exists $\epsilon>0$ small enough such that $\displaystyle{\gamma>\frac{b(K+\epsilon)}{1+E(K+\epsilon)}}.$ Hence, 
	$${\normalsize
	\begin{split}
		v'&=-v\left[\frac{\gamma+\delta v}{1+v}-b \frac{u}{1+Eu}\right] 
		\leq -v \left[\gamma- \frac{bu}{1+Eu}\right]
		\leq -v\left[\gamma- \frac{b(K+\epsilon)}{1+E(K+\epsilon)}\right] ,		
\end{split}}$$	
provided that $u(t)\leq K+\epsilon.$ Note that $u'<0$ if $u=K+\epsilon$ and $v\geq 0.$ Therefore, the set $\{(u,v)\in \mathbb{R}^{2}: 0<u\leq K+\epsilon, \;v>0\}$ is positively invariant. So, if the initial values satisfy $u(0)\leq K+\epsilon,$ $v(0)>0$ then $v\to 0$ exponentially as $t \to 0.$ On the other hand, if $u(0)>K+\epsilon$ then 
	$${\normalsize
	\begin{split}
		u'&=Au\left(1-\frac{u}{K}\right) - a\frac{uv}{1+Eu} 
		\leq Au
		\left(1-\frac{K+\epsilon}{K}\right)
		= -\frac{A \epsilon}{K}u ,		
\end{split}}$$	
provided that $u(t)>K+\epsilon.$ Hence, $u$ will be equal to $K+\epsilon$ in finite time and again $v\to 0$ as $t\to 0.$ Following a similar argument as the one in Lemma \eqref{asympt} we complete the proof in this case.

Next, we suppose that
$$(b-\gamma E)K=\gamma \Longleftrightarrow \gamma=\frac{b}{E}\frac{K}{\left(\frac{1}{E}+K\right)}.$$
Substituting this value of $\gamma$ in \eqref{KineHollinhII} and moving the origin into $(K,0)$ by the coordinate transformation $\tilde{u}=u-K,$ $\tilde{v}=v,$ one can rewrite the system \eqref{KineHollinhII} in the following form
  \begin{equation}\label{lia}
	\left \{
	\begin{array}{llll}
		\displaystyle{\tilde{u}' =  -A(\tilde{u}+K)\frac{\tilde{u}}{K}- \frac{a}{E} \frac{(\tilde{u}+K)\tilde{v}}{\left(\frac{1}{E} +\tilde{u}+K\right)}},\\
		\\
		\displaystyle{	\tilde{v}' = -\left[  \frac{\frac{b}{E}\frac{K}{\left(\frac{1}{E}+K\right)} + \delta \tilde{v}}{1+\tilde{v}}\right]\tilde{v} +\frac{b}{E}\frac{(\tilde{u}+K)\tilde{v}} {\left(\frac{1}{E}+\tilde{u}+K\right)}}.\\
	\end{array}
	\right.
\end{equation}
Now, we use the positive definite Liapunov function
$\displaystyle{V(\tilde{u}, \tilde{v})=\frac{b}{EK}\tilde{u}^{2}+ \left(\frac{1}{E} +K\right) \tilde{v}^{2}}.$
Denoting by $V'$ the derivative of $V$ with respect to the system \eqref{lia}, we obtain that
$${\small
\begin{split}
	-\frac{1}{2} V'(\tilde{u}, \tilde{v}) \left(\frac{1}{E}+\tilde{u}+K\right)(1+\tilde{v})&= \tilde{u}^{2}(\tilde{u}+K) \left(\frac{1}{E}+\tilde{u}+K\right)
	(1+\tilde{v})\frac{bA}{EK^{2}} 
	+
	\tilde{u} \tilde{v}(\tilde{u}+K)(\tilde{v}+1) \frac{ab}{E^{2}K} \\
	&\quad
	+\tilde{v}^{2}\left[\frac{bK}{E}+ \delta \tilde{v} \left(\frac{1}{E}+K\right)\right]
	\left(\frac{1}{E}+\tilde{u}+K\right)\\
	&
	\quad-\tilde{v}^{2}\frac{b}{E}(\tilde{u} +K)
	(1+\tilde{v})\left(\frac{1}{E}+K\right).
\end{split}
}
$$
Using that $\displaystyle{\delta \geq \frac{b}{E}}$ and performing a simple computation one  infers that
$${\small
	\begin{split}
		-\frac{1}{2} V'(\tilde{u}, \tilde{v}) \left(\frac{1}{E}+\tilde{u}+K\right)(1+\tilde{v})&\geq \tilde{u}^{2}(\tilde{u}+K) \left(\frac{1}{E}+\tilde{u}+K\right)
		(1+\tilde{v})\frac{bA}{EK^{2}} 
		+
		\tilde{u} \tilde{v}(\tilde{u}+K)(\tilde{v}+1) \frac{ab}{E^{2}K} \\
		&\quad
		+\frac{b}{E^{2}}\tilde{v}^{3} \left(\frac{1}{E}+K\right)-\frac{b}{E^{2}} \tilde{v}^{2} \tilde{u}.
	\end{split}
}
$$
Since $a\geq 1,$ we have that $V'(\tilde{u}, \tilde{v})<0,$ for $\tilde{u} \geq 0,$ and $\tilde{v}>0.$ Therefore, all solutions with positive initial conditions either tend to $(\tilde{u}, \tilde{v})=(0,0)$ or leave the $\tilde{u}\geq 0,$ $\tilde{v}>0$ quadrant through the line $\tilde{u}=0$ in finite time. On the other hand, the strip $\displaystyle{\{(\tilde{u},\tilde{v}):-K<\tilde{u}<0, \;p>0\}}$ is positively invariant and if $-K<\tilde{u}(t)<0,$ then
$${\normalsize
	\begin{split}
		\tilde{v}'(t) &\leq  -\frac{b}{E}\frac{\tilde{v}(t)}{(1+ \tilde{v}(t))} \left[  \frac{K}{\left(\frac{1}{E}+K\right)} -\frac{(\tilde{u}(t)+K)} {\left(\frac{1}{E}+\tilde{u}(t)+K\right)}+ \tilde{v}(t) -\frac{(\tilde{u}(t)+K)\tilde{v}(t)} {\left(\frac{1}{E}+\tilde{u}(t)+K\right)}\right]<0.
\end{split}
}
$$
Hence, once in the strip, $\tilde{v}(t)$ is monotone decreasing and $\tilde{v}(t) \to \beta \geq 0,$ as $t\to +\infty.$ Note that $\beta=0.$ In fact, if $\beta>0,$ then 
$${\normalsize
	\begin{split}
		\tilde{v}'(t) &< -\frac{b}{E}\frac{\tilde{v}(t)}{(1+ \tilde{v}(t))} \left[  1 -\frac{(\tilde{u}(t)+K)} {\left(\frac{1}{E}+\tilde{u}(t)+K\right)}\right]\tilde{v}(t)
		< -\frac{b}{E^{2}}\frac{\beta}{(1+ \tilde{v}(t_{0}))}   \frac{1} {\left(\frac{1}{E}+K\right)} \tilde{v}(t)		
	\end{split}
}
$$
would hold for some $t_{0}>0,$ and this would imply that $\tilde{v}\to 0$ exponentially, contradicting the assumption $\beta>0.$ Therefore,  $\tilde{v}(t) \to 0,$ as $t\to +\infty$ and the argument of the previous Lemma can be repeated again. This proves the Lemma.
\end{proof}
Depending on the parameters, the system \eqref{KineHollinhII} has at least one equilibrium with positive coordinates which can be found if we rewrite \eqref{KineHollinhII} as
  \begin{equation}\label{KineHollinhII2}
	\left \{
	\begin{array}{llll}
		\displaystyle{u' = \frac{a}{b}h(u)[f(u)-v]},\\
		\\
		\displaystyle{	v' = v\left[ h(u)-\left(\delta + \frac{\gamma - \delta}{1+v}\right)\right]},\\
	\end{array}
	\right.
\end{equation}
where
$\displaystyle{h(u)=b\frac{u}{1+Eu}}$ and $\displaystyle{f(u)= \frac{A}{aK}(K-u)(1+Eu).}$ 
Therefore, the (nontrivial) critical points $(u,v)$ are obtained as the intersection of the prey null-cline 
\begin{equation}\label{preynulclin}
	\mathcal{C}_{1}:v=f(u)=\frac{A}{aK}[-Eu^{2}+(KE-1)u+K],
\end{equation}
 and the predator null-cline
 \begin{equation}\label{predanulcline}
 	\mathcal{C}_{2}:v=M^{-1}(h(u))=-c \frac{u-d}{u-e}=-c+\frac{c(e-d)}{e-u},
 	\end{equation}
where
$$c=\frac{b-\gamma E}{b- \delta E},\;d=\frac{\gamma }{b- \gamma E},\;e=\frac{\delta}{b-\delta E}.
$$
and $M$ is defined as in \eqref{mort}.
It is easy to show that the Jacobian matrix $A$ associated with \eqref{KineHollinhII2} at any equilibrium point $(u,v)\in \mathcal{C}_{1} \cap \mathcal{C}_{2} \cap {R}^{2}_{+}$ is given by
\begin{equation}\label{marix}
 A= \left(
\begin{array}{cc}
	\displaystyle{	\frac{a}{b}\left[h(u)f'(u)\right]}& \displaystyle{-\frac{a}{b}h(u)}\\
	v h'(u) &-vM'(v)\\
\end{array}
\right),
\end{equation}
where 
$\displaystyle{h'(u)=\frac{b}{(1+Eu)^{2}}},$ $\displaystyle{f'(u)=\frac{A}{aK}(KE-1-2Eu),}$ and $\displaystyle{M'(v)= -\frac{\gamma-\delta}{(1+v)^{2}}.}$
 
 Next, we analyze the existence and stability of the nontrivial equilibria $(u_{0},v_{0})\in \mathbb{R}^{2}_{+}$ of system \eqref{KineHollinhII} given by the intersection of the curves \eqref{preynulclin} and \eqref{predanulcline}. From \eqref{preynulclin} we infer that any nontrivial equilibrium has to satisfy the condition $0<u<K.$ Since we are interested in the case $0<\gamma<\delta,$ we obtain $b-\delta E< b-\gamma E.$ Thus, we have the following three cases:\\
 
 \noindent
 {\bf{Case I.}}
 $0<b-E\delta<b-E\gamma<b.$ Therefore, $c>1$ and $e>d>0.$ \\
If $d\geq K,$ then do not exists non-trivial equilibria important from a Biological point of view and $(K,0)$ is local  asymptotically stable (global under additional conditions see Lemma \ref{asympt}).
If $d<K$ then $(K,0)$ is a saddle point and there exists just one nontrivial equilibrium
 \begin{equation}\label{equo}
(u_{0},v_{0})=\left(u_{0}, \frac{A}{aK}(K-u_{0})(1+Eu_{0})\right) = \left(u_{0},-c\frac{(u_{0}-d)}{(u_{0}-e)}\right)\in  \mathcal{C}_{1} \cap \mathcal{C}_{2} \cap \mathbb{R}^{2}_{+}.
\end{equation}
Introducing the variable $u_{0}$ we see that the dependence of that equilibria on the parameters is best characterized by identifying this equilibria with the points of a 8-dimensional manifold in $(A, K, a, E, \delta, \gamma, b, u_{0})-$space defined by the equation 
$$
\mathcal{S}:\frac{A}{aK}(K-u_{0})(1+Eu_{0}) =-c\frac{(u_{0}-d)}{(u_{0}-e)},
$$
where $0<\gamma< \delta,$ $0<b-E\delta<b-E\gamma,$ $d<K$ and $0<u_{0}<K.$
The Jacobian matrix given in \eqref{marix} at the equilibrium $(u_{0},v_{0})$ takes the form
 \begin{equation}\label{marix2}
  {\scriptsize
 	\begin{split}
 A_{u_{0},v_{0}}&=\left(
 \begin{array}{cc}
 	\displaystyle{\frac{Au_{0}2E}{K(1+Eu_{0})} \left[\frac{KE-1}{2E}-u_{0}\right]}& \displaystyle{-\frac{au_{0}}{(1+Eu_{0})}}\\
 	\displaystyle{v_0\frac{b}{(1+Eu_{0})^{2}}} & \displaystyle{-v_{0}\frac{(\delta-\gamma)}{\left(1+v_{0}\right)^{2}}}\\
 \end{array}
 \right) = \left(
 \begin{array}{cc}
 	\displaystyle{\frac{Au_{0}2E}{K(1+Eu_{0})} \left[\frac{KE-1}{2E}-u_{0}\right]}& -\\
 	+ & -\\
 \end{array}
 \right)
\end{split}}
\end{equation}
and from the convexity of the function \eqref{predanulcline} we infer $d<u_{0}<e.$ Note that if
 $u_{0}\geq\frac{KE-1}{2E},$ then it is easy to show that the characteristic polynomial
	\begin{equation}\label{charact1} P(\rho)=\rho^{2}-\text{Tr}(A_{u_{0},v_{0}})\rho+\text{det}(A_{u_{0},v_{0}})
	\end{equation}
 associated to \eqref{marix2} satisfy $\text{Tr}(A_{u_{0},v_{0}})<0$ and $\text{det}(A_{u_{0},v_{0}})>0.$ Therefore, $P(\rho)$ is a Hurwitz's  polynomial. Clearly, in this case $(u_{0},v_{0})$ given by \eqref{equo} is local asymptotically stable and the matrix $A_{u_{0},v_{0}}$ is not excitable (see \eqref{excita2x2}). \\
 
\begin{figure}[h!]
	\begin{center}
		\includegraphics[width=6.5cm]{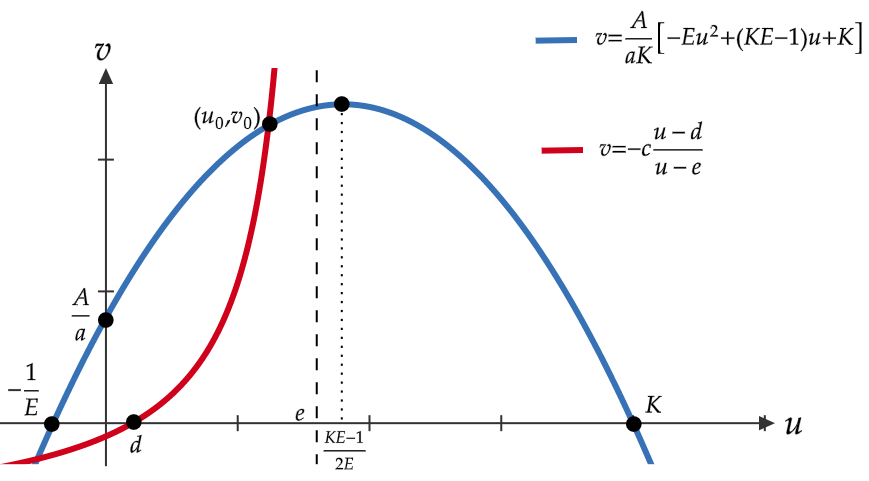}
	\caption{Location of equilibria if $\footnotesize{0<u_{0}< \frac{KE-1}{2E}.}$}\label{figura6}
\end{center}
\end{figure}

On the other hand, if
$0<u_{0}< \frac{KE-1}{2E},$ then the intersection point $(u_{0},v_{0})$ given by \eqref{equo} (see Figure \ref{figura6}) is local asymptotically stable if  \eqref{charact1} satisfies
\begin{equation}\label{tra1}
	\begin{split}
	\text{Tr}(A_{u_{0},v_{0}})<0& \Longleftrightarrow \frac{A}{EK} \frac{Eu_{0}}{(1+Eu_{0})}\left[E(K-u_{0})-(1+Eu_{0})\right]-\frac{v_{0}(\delta-\gamma)}{(1+v_{0})^{2}}<0,
	\end{split}
\end{equation}  
and
\begin{equation}\label{det1}
	\begin{split}
		\text{det}(A_{u_{0},v_{0}})>0& \Longleftrightarrow -\frac{Au_{0}}{K(1+Eu_{0})} \left[E(K-u_{0})-(1+Eu_{0}) \right]\frac{v_{0}(\delta-\gamma)}{(1+v_{0})^{2}}+\frac{abu_{0}v_{0}}{(1+Eu_{0})^{3}}>0.
	\end{split}
\end{equation} 
This last case is more interesting for our purpose  since conditions \eqref{tra1} and \eqref{det1} imply the matrix $A_{u_{0},v_{0}},$ given by  \eqref{marix2}, to be excitable (see \eqref{excita2x2}).
\begin{rem}
	Note that condition $KE-1>0$ ensures that for the prey an All\'ee effect zone exists, where the increase of prey density is favorable to its growth rate.
\end{rem}

From \eqref{tra1}--\eqref{det1} we infer that if the conversion rate $b$ and the gap between the maximal mortality $\delta$ and the mortality at low density $\gamma$ are big enough then the equilibrium $(u_{0},v_{0})$ given by \eqref{equo} is always local asymptotically stable. That is established in the following Theorem.  
\begin{thm}\label{allwsta}
Let the positive constants $A, K, a,$  and $E$ of system \eqref{KineHollinhII} be given with  $KE-1>0.$ Assume that $\gamma, \delta,$ and $b$ are given in such a way that $0<\gamma<\delta$ and $0<b-E\delta<b- E\gamma.$  If
	\begin{itemize}
		\item [$i)$]$\displaystyle{\delta-\gamma\geq a \left[1+\frac{A(KE+1)}{2a}\right]^{2}}$ 
		\item [$ii)$] $b>\displaystyle{\max\left\{\frac{\delta E(1+K)}{KE-1}, \frac{E^{2}(\delta-\gamma)(KE+1)}{8\left(E+\frac{A}{aK}\right)}\right\}}$ 
	\end{itemize}
	then the matrix $A_{u_{0},v_{0}},$ given by  \eqref{marix2}, is excitable (see \eqref{excita2x2}) and the unique nontrivial equilibrium 
	$(u_{0},v_{0}) \in  \mathcal{C}_{1} \cap \mathcal{C}_{2} \subset \mathbb{R}^{2}_{+}$ of system \eqref{KineHollinhII} is local asymptotically stable.
\end{thm}
\begin{proof}
	From hypothesis $b> \frac{\delta E(1+K)}{KE-1}.$ Therefore  $e<\frac{KE-1}{2E},$ and we infer $$0<d<u_{0}<e<\frac{KE-1}{2E}<K.$$ On the other hand, since $ u_{0}<\frac{KE-1}{2E}$	then it is easy to show that 
	\begin{equation}\label{vest}
		v_{0}=\frac{A}{a}\frac{(K-u_{0})}{K} (1+Eu_{0})<\frac{A}{a}\frac{(KE+1)}{2}.
	\end{equation}
Note that
$${\normalsize
\begin{split}
	\text{Tr}(A_{u_{0},v_{0}})&=\frac{A}{aK} \frac{Eu_{0}}{(1+Eu_{0})}(K-u_{0})(1+Eu_{0}) \frac{a}{(1+Eu_{0})}-\frac{A}{K}u_{0}
	-\frac{v_{0}(\delta-\gamma)}{(1+v_{0})^{2}}\\
	&< av_{0}\frac{Eu_{0}}{(1+Eu_{0})} \frac{1}{(1+Eu_{0})}-\frac{v_{0}(\delta-\gamma)}{(1+v_{0})^{2}}\\
	&< av_{0}-\frac{v_{0}(\delta-\gamma)}{\left(1+v_{0}\right)^{2}}.
\end{split}}$$
Therefore,
\begin{equation}\label{vest1}
	{\footnotesize
	\begin{split}
		\text{Tr}(A_{u_{0},v_{0}})<0& \Longleftarrow a-\frac{(\delta-\gamma)}{\left(1+v_{0}\right)^{2}}<0
\Longleftrightarrow 
v_{0}^{2}+2v_{0}+1-\frac{(\delta-\gamma)}{a}<0
 \Longleftrightarrow 
v_{0}<-1+\sqrt{\frac{\delta-\gamma}{a}}.
\end{split}}
\end{equation}
From \eqref{vest}  and \eqref{vest1} we deduce that a sufficient condition to obtain $\text{Tr}(A_{u_{0},v_{0}})<0$ is 
$$\frac{A}{a}\frac{(KE+1)}{2}\leq -1+\sqrt{\frac{\delta-\gamma}{a}}.$$
The last inequality holds from condition $i)$ of the Theorem. This implies \eqref{tra1}. Finally, using \eqref{equo}, we infer
$${\normalsize
\begin{split}
	\text{det}(A_{u_{0},v_{0}})
	&=\frac{au_{0}v_{0}}{(1+Eu_{0})^{2}}\left[-\frac{E(\delta-\gamma)v_{0}}{(1+v_{0})^{2}}	+\frac{A(\delta-\gamma)(1+Eu_{0})^{2}}{aK(1+v_{0})^{2}}
	+\frac{b}{(1+Eu_{0})}\right].\\
\end{split}
}$$
Since $\displaystyle{0<u_{0}<\frac{KE-1}{2E}},$ then we have
$${\normalsize
\begin{split}
	-\frac{E(\delta-\gamma)v_{0}}{(1+v_{0})^{2}}	+\frac{A(\delta-\gamma) (1+Eu_{0})^{2}}{aK(1+v_{0})^{2}}
	+\frac{b}{(1+Eu_{0})}
	&>	-\frac{E(\delta-\gamma)v_{0}}{(1+v_{0})^{2}}
	+\frac{A(\delta-\gamma)}{aK(1+v_{0})^{2}}
	+\frac{2b}{(1+EK)}.
\end{split}
}
$$
Thus, 
\begin{equation}\label{vest2}
	{\small
	\begin{split}
		\text{det}(A_{u_{0},v_{0}})>0& \Longleftarrow 	-\frac{E(\delta-\gamma)v_{0}}{(1+v_{0})^{2}}
		+\frac{A(\delta-\gamma)}{aK(1+v_{0})^{2}}	
		+\frac{2b}{(1+EK)}>0	\\
		& \Longleftrightarrow 
	\frac{2b}{(KE+1)} v_{0}^{2}+ \left(\frac{4b}{(KE+1)}-E(\delta-\gamma)\right)v_{0} 
	+\frac{A(\delta-\gamma)}{aK}
	+\frac{2b}{(KE+1)}>0.
	\end{split}}
\end{equation}
Condition $ii)$ of the Theorem assures 
the quadratic polynomial in the variable $v_{0}$ given by  \eqref{vest2} has complex roots and \eqref{det1} holds.
\end{proof}

\begin{figure}[h!]
 		\begin{center}
 			\includegraphics[width=7cm]{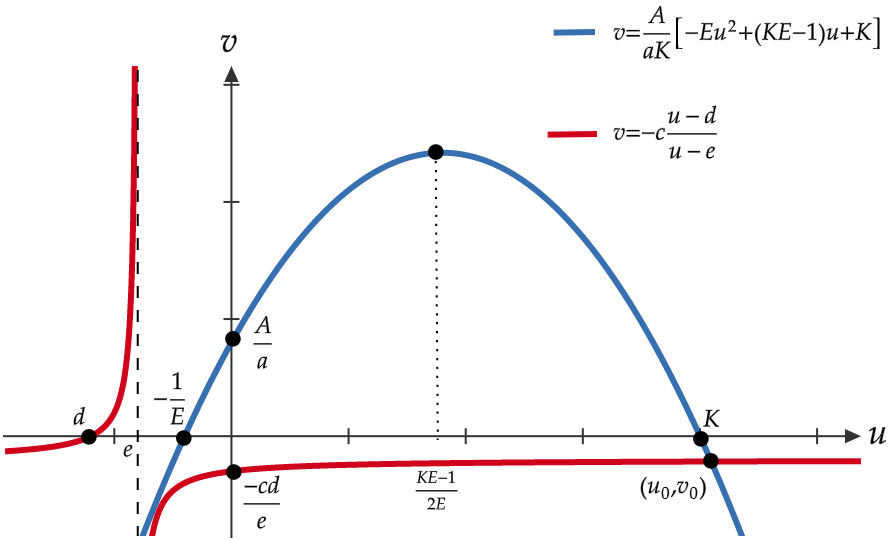}
 			\caption{Location of equilibria under assumptions $b-\delta E< b- \gamma E<0.$  }\label{figura7}
 		\end{center}
 	\end{figure}

\noindent
{\bf{Case II.}} 	 $b-\delta E< b- \gamma E<0.$ Thus \eqref{lasympsta} implies that $(K,0)$ is local (global under additional conditions see Lemma \eqref{asympt}) asymptotically stable. In this case $0<c<1,$ $d<e<-\frac{1}{E}<0.$ Therefore the system \eqref{KineHollinhII} does not have a  relevant nontrivial equilibrium from a biological point of view (see Figure \ref{figura7}).\\
 	
 \noindent
 {\bf{Case III.}} $b-E\delta<0<b-E\gamma<b.$
 This condition implies that $c<0,$ $e<-\frac{1}{E}<0<d.$ \\
 If $d\geq K,$ the equilibria $(K,0)$ is local  asymptotically stable (global under additional conditions see Lemma \ref{globsta}) and the system \eqref{KineHollinhII} does not have a  relevant nontrivial equilibrium  from a biological point of view.\\
 
 If $d< K$  then $(K,0)$ is a saddle point and we will show that under a suitable choice of parameters the system \eqref{KineHollinhII} can have one, two or three nontrivial equilibria. 
 In fact, It is easy to show that $$f(u^{\ast})=\frac{AE}{aK}\left(K-u^{\ast}\right) \left(\frac{1}{E}+u^{\ast}\right)=
  \frac{AE}{4aK} \left(\frac{1}{E}+K\right)^{2},$$
 where 
 $u^{\ast}= \displaystyle{\frac{KE-1}{2E}}$ is the abscissa of the parabola's vertex defined by \eqref{preynulclin}. We study two sub-cases:
 
  \noindent
 {\bf{$i)$}} $\displaystyle{KE-1\leq 0}.$
 	Then, there exists just one nontrivial local asymptotically stable equilibrium important from a Biological point of view 
 	$(u_{0},v_{0})\in  \mathcal{C}_{1} \cap \mathcal{C}_{2} \cap \mathbb{R}^{2}_{+}.$
In this case $d<u_{0}<K,$ the Jacobian matrix $A_{u_{0},v_{0}},$ given by \eqref{marix2}, satisfy relations \eqref{tra1}-\eqref{det1} and it is not excitable.

\begin{figure}[h!]
	\begin{center}
		\includegraphics[width=8cm]{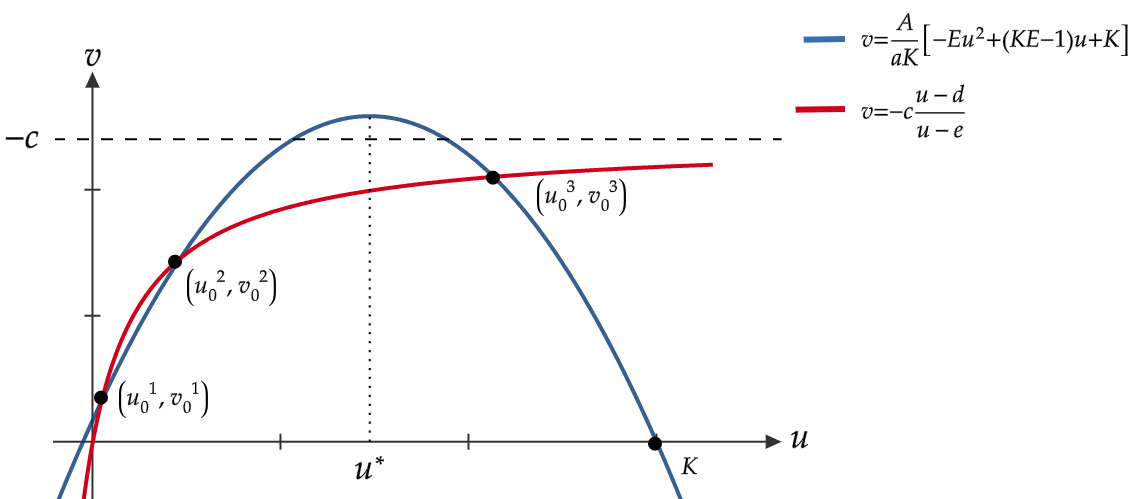}
		\caption{Location of equilibria under assumptions $b-\delta E<0< b- \gamma E<0.$  }\label{figura8}
	\end{center}
\end{figure}

 \noindent
{\bf{$ii)$}} $KE-1> 0.$
Let us choose the parameters of system 
\eqref{KineHollinhII2} in such a way that
\begin{equation}\label{ord}
	2E\gamma <b, \quad  KE>3,
\end{equation}
\begin{equation}\label{ord1}
f\left(\frac{1}{E}\right)<M^{-1}\left(h \left(\frac{1}{E}\right)\right),
\end{equation}
 and
 \begin{equation}\label{ord2}
 \frac{AKE}{4a}  >-c.
 \end{equation}
From \eqref{ord}, we infer that 
$$d<\frac{1}{E}<u^{\ast}.$$ 
Therefore, \eqref{ord1} implies that there exists $u_{0,1} \in \left(d, \displaystyle{\frac{1}{E}}\right)$ such that $f\left(u_{0,1}\right)=M^{-1}\left(h \left( u_{0,1}\right)\right).$ On the other hand, from \eqref{ord2} we have that
$$f(u^{\ast})  > \frac{AKE}{4a}  >-c.$$ 
  Hence, it follows that there exists $u_{0,2} \in \left(\frac{1}{E}, u^{\ast}-\frac{\sqrt{2KE +1}}{2E}\right)$ and $u_{0,3} \in (u^{\ast}, K)$ such that $f\left(u_{0,i}\right)=M^{-1}\left(h \left( u_{0,i}\right)\right),$ $i=2,3.$ Therefore, $(u_{0,i}, f(u_{0,i}))=(u_{0,i}, M^{-1}(h(u_{0,i}))) \in \mathcal{C}_{1} \cap \mathcal{C}_{2} \cap \mathbb{R}^{2}_{+},$ $i=1,2,3$ (see \eqref{preynulclin}-\eqref{predanulcline}) are three nontrivial equilibria  (see Figure \ref{figura8}).
From this analysis, it is not difficult to infer that $u_{0,1}$ and $u_{0,2}$ ($u_{0,2}$ and $u_{0,3}$) can collapse generating two nontrivial equilibrium points. This proves our assertion. Regarding their stability, similar to the part $i)$ of {\bf{Case I,}} we can show that $(u_{0,3},v_{0,3})$ is local asymptotically stable but the matrix $A_{u_{0,3},v_{0,3}}$ is not excitable (see \eqref{excita2x2}). On the other hand, unlike {\bf{Case I,}} it is difficult to prove the local stability of equilibrium $(u_{0,i},v_{0,i}),$ for $i=1,2.$ Nevertheless, some some simulations has shown that $(u_{0,1},v_{0,1})$ is a stable node  and $(u_{0,2},v_{0,2})$ is a saddle point. This suggest that the matrix $A_{u_{0,1},v_{0,1}}$ defined similarly as in \eqref{marix2} has the possibility of being excitable (see \eqref{excita2x2}) and, therefore, the equilibrium $(u_{0,1},v_{0,1})$ could undergoes a Turing bifurcation. 

One can give sufficient conditions on the parameters of the original system \eqref{KineHollinhII} to assure inequalities \eqref{ord}-\eqref{ord1}-\eqref{ord2} hold. That is stated in the following Lemma.
\begin{lem}\label{para}
	Suppose that positive parameters $A, K, E, a, b, E, \gamma,$ and $\delta$ of system \eqref{KineHollinhII} satisfy:
  
  \begin{equation}\label{param1}
  0<\frac{A}{a}<\frac{1}{2},
  \end{equation}

   \begin{equation}\label{param21}
 KE>
 \max \left\{8, \frac{\left(8-\frac{A}{a}\right) + \sqrt{\frac{A^{2}}{a^{2}} + \frac{16A}{a} +48 } }{2 \left(\frac{1}{2}-\frac{A}{a}\right)} \right\},
 \end{equation}
 
 \begin{equation}\label{param13}
0<\frac{\gamma}{\delta}<\min\left\{ \frac{1}{2},\; \frac{1}{2}-\frac{2A u^{\ast}}{a K},\; \frac{AEK}{2(AEK+2a)},\; \frac{\frac{A}{a}\left[ \left(\frac{1}{2}-\frac{A}{a}\right)(KE)^{2} +  \left(\frac{A}{a}-8\right) KE+8\right]}{\frac{A}{a}(KE)^{2}+4 \left(\frac{1}{2}-\frac{A}{a}\right)KE+ \frac{4A}{a}} \right\},
\end{equation}
and
\begin{equation}\label{param2}
\max\left\{2E\gamma, \;2E\delta - \frac{2E(\delta-\gamma)}{\frac{4A u^{\ast}}{a K} + 1}\right\} <b < \min \left\{ E \delta,\; \frac{E\left(4a\gamma+AEK\delta \right)}{4a+AEK} \right\}.
\end{equation} 
Then conditions  \eqref{ord}, \eqref{ord1}, and \eqref{ord2} hold.
\end{lem}
\begin{proof}
Note that \eqref{param21} and \eqref{param1} imply $KE>8.$ From this and \eqref{param2} we obtain \eqref{ord}.  Also, from \eqref{param2}, we obtain \eqref{ord1}  (see \eqref{preynulclin}- \eqref{predanulcline}).
 In fact, 
 $${\normalsize
 	\begin{split}
	f\left(\frac{1}{E}\right)<M^{-1}\left(h \left(\frac{1}{E}\right)\right) &\Longleftrightarrow
\frac{4A u^{\ast}}{a K}< \frac{(b-2E \gamma)}{-(b-2E \delta)}\\
 & \Longleftrightarrow 
b>2E\delta - \frac{2E(\delta-\gamma)}{ \frac{4A u^{\ast}}{aK} +1}. 
 	\end{split}
}
 $$
Therefore, the last term in the above inequality must satisfy
  $${\normalsize
 	\begin{split}
 2E\delta - \frac{2E(\delta-\gamma)}{ \frac{4A u^{\ast}}{aK} +1}<E \delta
  &\Longleftrightarrow
 		\frac{\gamma}{\delta}<\frac{1}{2} -\frac{2A u^{\ast}}{a K}.
 	\end{split}
 }
 $$
 The last inequality holds because  \eqref{param13}. Note that \eqref{param1} assure that, in the last inequality, $\displaystyle{\frac{\gamma}{\delta}}$ is positive. On the other hand, from \eqref{param2} we have that
   $$\frac{AKE}{4a}  >-c=\frac{b-\gamma E}{-(b- \delta E)}.$$ 
   
   In fact,
 $${\normalsize
	\begin{split}
\frac{AKE}{4a}  >\frac{b-\gamma E}{-(b- \delta E)} &\Longleftrightarrow
\frac{AKE}{4a} > -1+ \frac{E(\delta -\gamma)}{-(b-E \delta)}\\
		& \Longleftrightarrow 
		b<\frac{E\left(4a\gamma+AEK\delta \right)}{4a+AEK}. 
	\end{split}
}
$$
 Therefore, the last term in the above inequality should satisfy
$${\normalsize
	\begin{split}
	\frac{E\left(4a\gamma+AEK\delta \right)}{4a+AEK}>2E \gamma
		&\Longleftrightarrow
		\frac{\gamma}{\delta}<\frac{AEK}{2(AEK+2a)},
	\end{split}
}
$$   
  and 
$${\normalsize
	\begin{split}
2E\delta - \frac{2E(\delta-\gamma)}{ \frac{4A u^{\ast}}{aK} +1}< \frac{E\left(4a\gamma+AEK\delta \right)}{4a+AEK}
&\Longleftrightarrow		
		\frac{\gamma}{\delta}< \frac{\frac{A}{a}\left[ \left(\frac{1}{2}-\frac{A}{a}\right)(KE)^{2} +  \left(\frac{A}{a}-8\right) KE+8\right]}{\frac{A}{a}(KE)^{2}+4 \left(\frac{1}{2}-\frac{A}{a}\right)KE+ \frac{4A}{a}}.
	\end{split}
}
$$   
\eqref{param13} ensures that both last inequalities in the right side last are satisfied and  using
\eqref{param1}--\eqref{param21} it is easy  to show that $\displaystyle{\frac{\gamma}{\delta}}>0.$   
  This proves the Lemma. 
 \end{proof}

 \subsection{Existence of patterns for the predator-prey model \eqref{FunctionResp}-\eqref{NewmanCond}.}\label{patt}
 Here we  study the effect of diffusion on the stability of equilibria in the reaction-diffusion model \eqref{FunctionResp}-\eqref{NewmanCond} and explore under which parameter values a Turing instability can occur. Furthermore, we apply the Criterion II (see Theorem \ref{Gensystem18}) to provide specific diffusive parameters values to ensure that a Turing bifurcation occurs for the system \eqref{FunctionResp}-\eqref{NewmanCond} giving rise to non-uniform stationary solutions.  That is established in the following two results.

 \begin{thm}\label{PatPP1}
 Let the positive parameters $A, K, a,E$ of system \eqref{FunctionResp}-\eqref{NewmanCond}  be given with  $KE-1>0.$  Assume the  parameters $\gamma,\delta$ and $b$ of system \eqref{FunctionResp}-\eqref{NewmanCond} satisfying:  
 \begin{itemize}
 	\item[$i)$] $0<E\gamma<E\delta<b,$
 	
 	\item[$ii)$] $\displaystyle{\delta-\gamma\geq a \left[1+\frac{A(KE+1)}{2a}\right]^{2}},$ 
 	and		
 	\item [$iii)$] $b>\displaystyle{\max\left\{\frac{\delta E(1+K)}{KE-1}, \frac{E^{2}(\delta-\gamma)(KE+1)}{8\left(E+\frac{A}{aK}\right)}\right\}}.$ 
 \end{itemize}
Consider the matrix
	$(a_{ij})=A_{u_{0},v_{0}}$ defined as in \eqref{marix2}
where $w_{0}=(u_{0},v_{0})$ is the unique equilibrium point of system \eqref{KineHollinhII2} belonging to the intersection  $\mathcal{C}_{1} \cap \mathcal{C}_{2} \cap \mathbb{R}^{2}_{+}$ (see \eqref{preynulclin}--\eqref{predanulcline}) and 
 denote by $\displaystyle{\lambda_{j},}$ $j=0,1,2,...$ the eigenvalues, with respective eigenfunctions $\phi_{j},$  of Laplacian operator $-\Delta$ on  $\Omega,$ with  no flux-boundary conditions. Suppose that $\displaystyle{\lambda_{l},}$ is a simple eigenvalue for some  $l\in \mathbb{N},$ with $l \geq 2.$ 
 
 If 
$	d_{1}'<d_{1}<d_{1}'',$ where $d_{1}'$ and $d_{1}''$ are given by Remark \eqref{crtII}, then   at
 	$\displaystyle{d_{2}^{\ast}=\frac{a_{22}\lambda_{l}d_{1} -\text{det}(A_{u_{0},v_{0}})}{\lambda_{l}(\lambda_{l}d_{1}-a_{11})}}$ the uniform steady-state solution $w_{0}$ of 
 	\eqref{FunctionResp}-\eqref{NewmanCond}, undergoes a Turing bifurcation.
 	Furthermore, $$w(x,s)=w_{0}+s\cdot v_{1l}\phi_{l}(x)+O(s)^{2}$$
 	is a one-parameter family of non-uniform stationary solutions of \eqref{FunctionResp}-\eqref{NewmanCond},  with $s\in (-\zeta, \zeta),$ for some small enough $\zeta$ and
 	 $$v_{1l}= \left(
 	 \begin{array}{c}
 	 	-\frac{(a_{22}-\lambda_{l}d_{2}^{\ast})}{a_{21}}	\\
 	 	1\\
 	 \end{array}
 	 \right).$$
Moreover,  if $\lambda_{1}$ is a simple eigenvalue and  $\frac{a_{11}}{\lambda_{2}}\leq d_{1}<\frac{a_{11}}{\lambda_{1}}$ then at $\displaystyle{d_{2}^{\ast}=\frac{a_{22}\lambda_{1}d_{1} -\text{det}(A_{u_{0},v_{0}})}{\lambda_{1}(\lambda_{1}d_{1}-a_{11})}},$ the uniform steady-state solution $w_{0}$ of 
\eqref{FunctionResp}-\eqref{NewmanCond} undergoes a Turing bifurcation.
 \end{thm}

\begin{proof}
From the analysis carried out in {\bf{Case I}} of Subsection  \ref{Applications} we know that $w_{0}=(u_{0},v_{0})$ is a uniform steady state solution of system 	\eqref{FunctionResp}-\eqref{NewmanCond}. On the other hand, 	from Theorem \ref{allwsta} we have that the Jacobian matrix $A_{u_{0},v_{0}}$ defined as in \eqref{marix2} is excitable with $a_{11}>0.$ The result follows as a direct consequence of Theorem \ref{Gensystem18}.
\end{proof}

We have perform some numeric simulations to find the face of the  patterns for the system \eqref{FunctionResp}-\eqref{NewmanCond} using the hypotheses of Theorem \ref{PatPP1}. For this, we assume that  the function $d_{2}:(-\zeta, \zeta) \to \mathbb{R}_{+}$ given  in \eqref{Gensystem20} have the form $d_{2}(s)=d_{2}^{\ast}-s,$ where    $\zeta\leq \min\{d_{2}^{\ast}-d_{2}', \;d_{2}'''-d_{2}^{\ast}\}.$ Note that  $d_{2}(0)=d_{2}^{\ast}$  and $d_{2}(s)$ defined in this way, could not be the same function given by Theorem \ref{Gensystem18} in \eqref{Gensystem20}. Nevertheless, considering $d_{2}$ as a linear function is useful for performing simulations.

 Next, we suppose that $\Omega$ is one-dimensional given by the interval $(0,100).$  Therefore, we are interested in non-uniform stationary solutions $u:\mathbb{R}^{+} \times [0,100] \to \mathbb{R}^{+}$ and $v:\mathbb{R}^{+} \times [0,100] \to \mathbb{R}^{+}$ of system \eqref{FunctionResp}-\eqref{NewmanCond} that satisfy the no-flux boundary condition 
$$\frac{\partial u(t,0)}{\partial x}=
\frac{\partial u(t,100)}{\partial x}=0=\frac{\partial v(t,100)}{\partial x}
=\frac{\partial v(t,0)}{\partial x}.$$
Taking into account the domain $\Omega,$ we can rewrite the boundary problem \eqref{Gensyst7} as 
$$
	\phi''_{j}(x)+ \lambda_{j}\phi_{j}(x)=0,\;\; x\in(0,100);\;\;\; \;\;\;\;\;\; 
	\phi'_{j}(0)= \phi'_{j}(100) =0.
$$
Thus the (simple) eigenvalues are given by $\displaystyle{\lambda_{j}=\left(\frac{j \pi}{100} \right)^{2}},$ $j=0,1,2,...$
with corresponding eigenfunctions $\displaystyle{\phi_{j}(x)=\cos\left(\frac{j\pi x}{100}\right)}.$ It is easy to verify that parameters $A=3,$ $K=15,$ $a=2,$ $E=1/2,$ $\delta=129,$  $\gamma=20,$  $b=159,$ $\lambda_{10}=0.0799$ and $d_{1}=3.2945$ satisfy the hypotheses of Theorem \ref{PatPP1}. Thus, at $d_{2}^{\ast}=5358,$ the uniform steady-state solution $w_{0}=(0.8226, 2.0009) \in \mathcal{C}_{1} \cap \mathcal{C}_{2} \cap \mathbb{R}^{2}_{+}$ of system \eqref{FunctionResp}-\eqref{NewmanCond} undergoes a Turing bifurcation. If we fix $s=-0.1$ and $d_{2}=5358.1$ then we obtain that 
\begin{equation}\label{appr1}
\left(\tilde{u}(x), \tilde{v}(x)\right)= \left(0.8326-0.3462 \cos\left(\frac{\pi x}{10}\right), 2.0109-0.1 \cos\left(\frac{\pi x}{10}\right)\right)
\end{equation}
 is an approximation of a non-uniform stationary solution of \eqref{FunctionResp}-\eqref{NewmanCond}. Using \eqref{appr1} as initial data we see the evolution in time and the convergence when of the Midpoint Method to an exact non-uniform stationary solution of system  \eqref{FunctionResp}-\eqref{NewmanCond}, see \cite[\href{https://tinyurl.com/444wuvpm}{\blue Interactive Simulation 1}]{Walker et all} (the reader would be able to change the $v-$predator or $u-$prey solution displayed and other characteristics in the simulation in the pane on the left side).

 \begin{figure}[h!]
 	\begin{center}
 		\includegraphics[width=6.5cm]{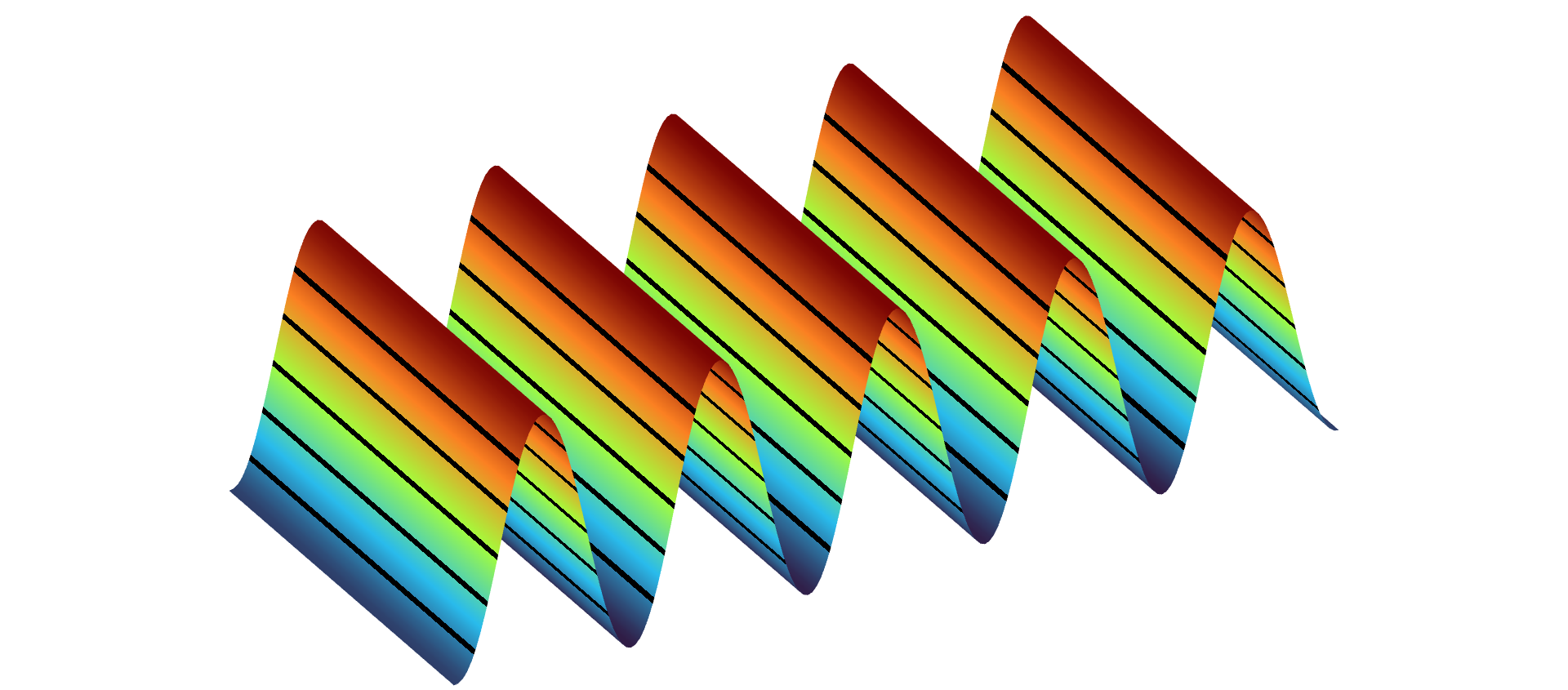}
 		\includegraphics[width=6.5cm]{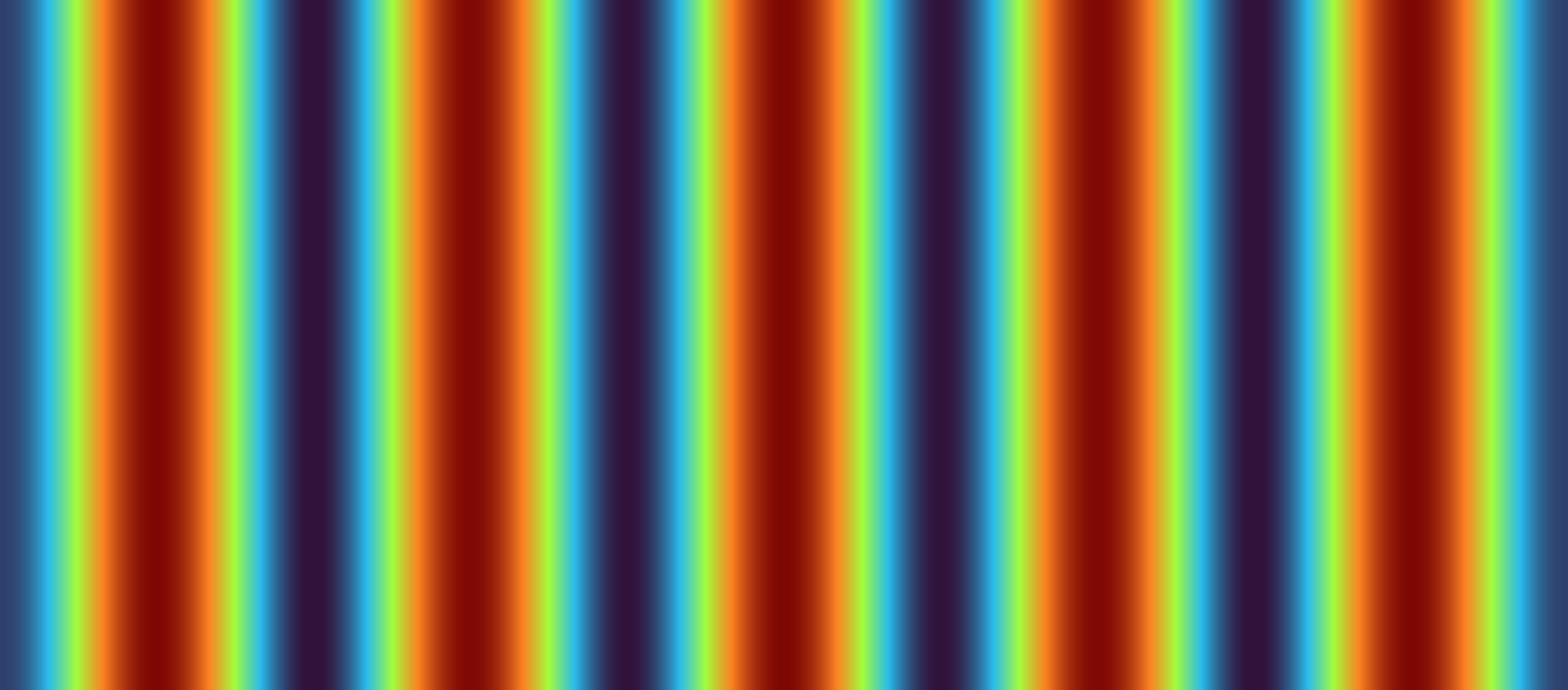}
 		\caption{The population density of predator $1.92 \leq v(x) \leq 2.12.$}\label{figura10}
 	\end{center}
 \end{figure}

 \begin{figure}[h!]
	\begin{center}
		\includegraphics[width=8cm]{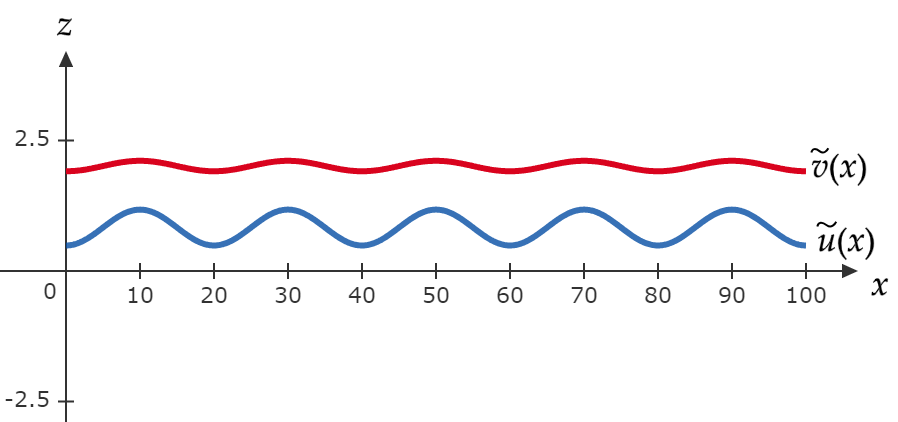}
		\caption{Projection of the population density approximations of predator $\tilde{v}$ and the prey $\tilde{u}$ on the $xz-$plane.}\label{figura11}
	\end{center}
\end{figure}

\begin{rem}
It is clear that conditions of Theorem \ref{PatPP1} are sufficient but not necessary, i.e., one can choose another appropriated distribution of parameters for system \eqref{FunctionResp}-\eqref{NewmanCond} that do not satisfy conditions  $ii)$ and $iii)$ of Theorem \ref{PatPP1} but still preserve the hypotheses of Criterion II to have the formation of patterns. For instance, choose the parameters as: $A=3,$ $K=15,$ $a=2,$ $E=1/2,$ $\gamma=9$ $\delta=78,$ $b=92,$ $\lambda_{10}=0.0799,$ and $d_1=3.3603.$ Hence at $d_{2}^{\ast}=3053.4$ the uniform steady-state solution $w_{0}=(0.8561, 2.0198) \in \mathcal{C}_{1} \cap \mathcal{C}_{2} \cap \mathbb{R}^{2}_{+}$ of system \eqref{FunctionResp}-\eqref{NewmanCond} undergoes a Turing bifurcation. If we fix $s=-0.1$ and $d_{2}=3053.5$ then we obtain that solution  
\begin{equation}\label{appr2}
	\left(\tilde{u}(x), \tilde{v}(x)\right)= \left(0.8661-0.3475 \cos\left(\frac{\pi x}{10}\right), 2.0298-0.1 \cos\left(\frac{\pi x}{10}\right)\right)
\end{equation}
is an approximation of a non-uniform stationary solution of system \eqref{FunctionResp}-\eqref{NewmanCond}, see \cite[\href{https://tinyurl.com/yzj8fsfv}{\blue Interactive Simulation 2}]{Walker et all} when $t\geq 55.$
\end{rem}
 
 \begin{thm}\label{PatPP2}
	Suppose that positive parameters $A, K, E, a, b, E, \gamma,$ and $\delta$ of system \eqref{FunctionResp}-\eqref{NewmanCond} satisfy conditions \eqref{param1}, \eqref{param21}, \eqref{param13} and \eqref{param2}.
	Consider the matrix
$(a_{ij})=A_{u_{0,1},v_{0,1}}$ defined as in \eqref{marix2} where $w_{0,1}=(u_{0,1},v_{0,1})$ is the equilibrium point of system \eqref{KineHollinhII2} belonging to the intersection  $\mathcal{C}_{1} \cap \mathcal{C}_{2} \cap \mathbb{R}^{2}_{+}$ (see \eqref{preynulclin}--\eqref{predanulcline}) with $u_{0,1} \in \left(d, \frac{1}{E}\right)$ and 
	denote by $\displaystyle{\lambda_{j},}$ $j=0,1,2,...$ the eigenvalues, with respective eigenfunctions $\phi_{j},$  of Laplacian operator $-\Delta$ on  $\Omega,$ with  no flux-boundary conditions. Suppose that $\displaystyle{\lambda_{l},}$ is a simple eigenvalue for some  $l\in \mathbb{N},$ with $l \geq 2.$
	If $\text{det}(A_{u_{0,1},v_{0,1}})>0,$ $\text{Tr}(A_{u_{0,1},v_{0,1}})<0,$ and
	$	d_{1}'<d_{1}<d_{1}'',$ where $d_{1}',$  $d_{1}''$ are given by Remark \eqref{crtII}, then   at
	$\displaystyle{d_{2}^{\ast}=\frac{a_{22}\lambda_{l}d_{1} -\text{det}(A_{u_{0,1},v_{0,1}})}{\lambda_{l}(\lambda_{l}d_{1}-a_{11})}}$ the uniform steady-state solution $w_{0,1}$ of system
	\eqref{FunctionResp}-\eqref{NewmanCond}, undergoes a Turing bifurcation.
	Furthermore, $$w(x,s)=w_{0,1}+s\cdot v_{1l}\phi_{l}(x)+O(s)^{2}$$
	is a one-parameter family of non-uniform stationary solutions of \eqref{FunctionResp}-\eqref{NewmanCond},  with $s\in (-\zeta, \zeta),$ for some small enough $\zeta$  and
	$$v_{1l}= \left(
	\begin{array}{c}
		-\frac{(a_{22}-\lambda_{l}d_{2}^{\ast})}{a_{21}}	\\
		1\\
	\end{array}
	\right).$$
	Moreover,  if $\lambda_{1}$ is a simple eigenvalue and  $\frac{a_{11}}{\lambda_{2}}\leq d_{1}<\frac{a_{11}}{\lambda_{1}}$ then and at $\displaystyle{d_{2}^{\ast}=\frac{a_{22}\lambda_{1}d_{1} -\text{det}(A_{u_{0,1},v_{0,1}})}{\lambda_{1}(\lambda_{1}d_{1}-a_{11})}},$ the uniform steady-state solution $w_{0,1}$ of system
	\eqref{FunctionResp}-\eqref{NewmanCond} undergoes a Turing bifurcation.
\end{thm}

\begin{proof} The result is consequence of the analysis carried out in {\bf{Case III}} of subsection  \ref{Applications} and similar arguments as in the proof of Theorem \ref{PatPP1}. 
\end{proof}
 
\begin{figure}[h!]
	\begin{center}
		\includegraphics[width=5cm]{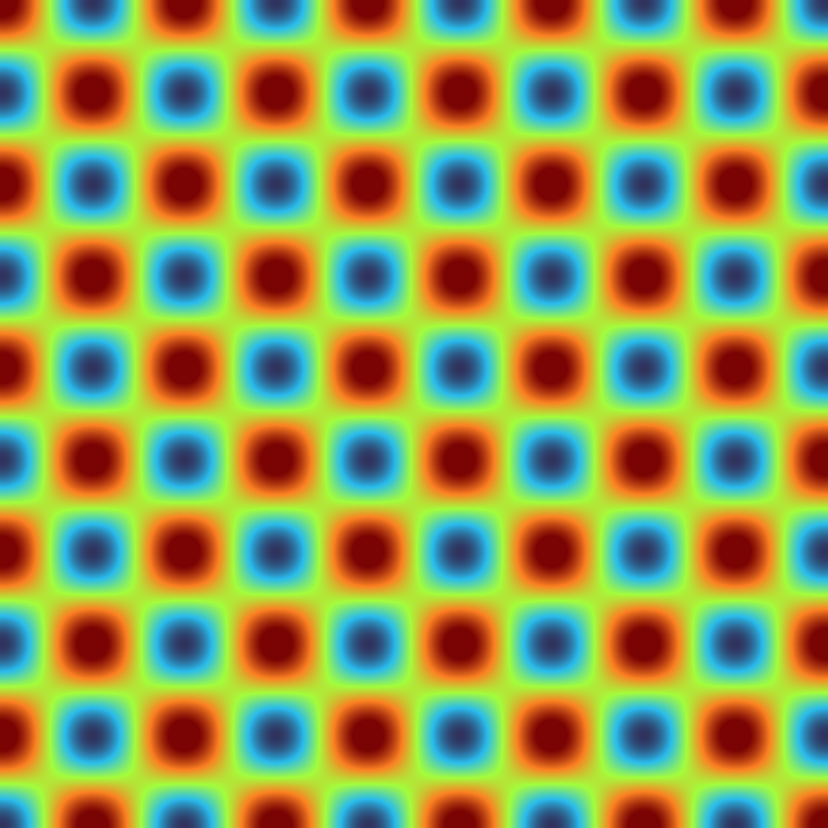}
		\includegraphics[width=6cm]{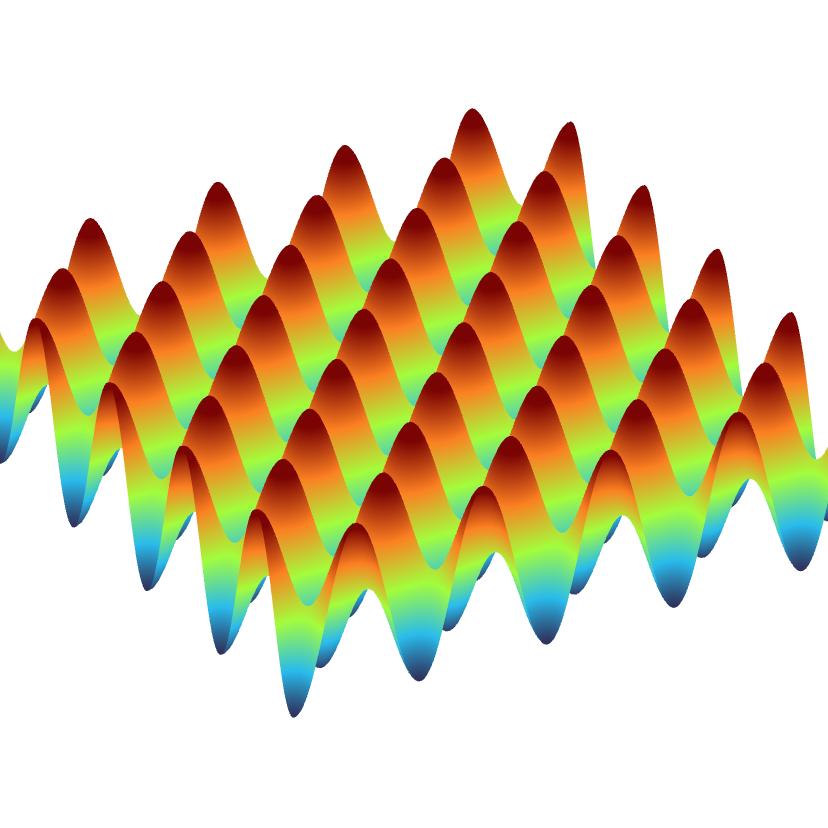}
		\caption{The population density of prey $0.65 \leq v(x) \leq 0.656.$}\label{figura15}
	\end{center}
\end{figure}

Next, we perform some simulations, considering hypotheses of Theorem \ref{PatPP2}. In this case, we consider the bidimensional domain $\Omega=(0,L_{1}) \times (0,L_{2}),$ with  $L_{1}>L_{2}$ and we suppose that  $\frac{L_{2}^{2}}{L_{1}^{2}}$ is an irrational number. Therefore, we are interested in solutions $u:\mathbb{R}^{+} \times [0,L_{1}]\times [0,L_{2}] \to \mathbb{R}^{+}$ and $v:\mathbb{R}^{+} \times [0,L_{1}]\times [0,L_{2}] \to \mathbb{R}^{+}$ of system \eqref{FunctionResp}-\eqref{NewmanCond} that satisfy the no-flux boundary conditions
$$\frac{\partial u(t,0,y)}{\partial x}=
\frac{\partial u(t,L_{1},y)}{\partial x}=0=\frac{\partial v(t,L_{1},y)}{\partial x}
=\frac{\partial v(t,0,y)}{\partial x},$$

$$\frac{\partial u(t,x,0)}{\partial y}=
\frac{\partial u(t,x,L_{2})}{\partial y}=0=\frac{\partial v(t,x,L_{2})}{\partial y}
=\frac{\partial v(t,x,0)}{\partial y}.$$
In this case, we can rewrite the boundary problem \eqref{Gensyst7} as 
\begin{equation}\label{unidlapl4}
	\begin{split}
		\frac{\partial^{2} \phi(x,y)}{\partial x^{2}}	+ \frac{\partial^{2} \phi(x,y)}{\partial y^{2}}+ \lambda\; \phi(x,y)=0,&\;\; x\in(0,L_{1})\times (0,L_{2});\\
		\frac{\partial \phi(0,y)}{\partial x}=0= \frac{\partial \phi(L_{1},y)}{\partial x},&\;\; 0\leq y \leq L_{2}.\\
		\frac{\partial \phi(x,L_{2})}{\partial y}=0=
		\frac{\partial \phi(x,0)}{\partial y},& \;\;0\leq x \leq L_{1}.\\
	\end{split}
\end{equation}
Hence, the eigenvalues are given by $\displaystyle{\tilde{\lambda}_{n,m}=\pi^{2} \left(\frac{n^{2}}{L_{1}^{2}}+\frac{m^{2}}{L_{2}^{2}}\right)},$ $n,m=0,1,2,...$
with corresponding eigenfunctions $\displaystyle{\phi_{n,m}(x,y)=\cos\left(\frac{n\pi x}{L_{1}}\right) \cos\left(\frac{m\pi y}{L_{2}}\right)}.$ Irreducibility of $\frac{L_{2}^{2}}{L_{1}^{2}}$ ensures that the sequence of indices $\{(n,m)\}_{(n,m)\in \mathbb{N}\times \mathbb{N}}$ can be arranged in order that $\lambda_{j}:=\tilde{\lambda}_{n_{j},m_{j}},$ $j=0,1,2,...$ is simple and $0=\lambda_{0}<\lambda_{1}=<\lambda_{2}<\cdots<\lambda_{j}<\cdots.$ Choosing $L_{1}=100 \sqrt[4]{2},$ $L_{2}=50,$ $A=1,$ $K=25,$ $a=3,$ $E=2,$ $\delta=100,$  $\gamma=1,$  $b=160.5,$ $\lambda_{99}:=\tilde{\lambda}_{8,9}\approx 0.3644$ and $d_{1}=0.4622$ it is easy to prove that hypotheses of Theorem \ref{PatPP2} are satisfied and therefore at $d_{2}^{\ast}=100.0809,$ the uniform steady-state solution $w_{0}=(0.4998, 0.6532) \in \mathcal{C}_{1} \cap \mathcal{C}_{2} \cap \mathbb{R}^{2}_{+}$ of system \eqref{FunctionResp}-\eqref{NewmanCond} undergoes a Turing bifurcation. If we fix $s=-0.02$ and $d_{2}=100.1009$ then we obtain that 
\begin{equation}\label{appr3}
	{\footnotesize
	\begin{split}
	\left(\tilde{u}(x), \tilde{v}(x)\right)&= \left(0.5002-0.0497 \cos\left(\frac{ 8 \pi x}{100 \sqrt[4]{2}}\right) \cos\left(\frac{9 \pi y}{50}\right), 0.6536-0.02 \cos\left(\frac{ 8 \pi x}{100 \sqrt[4]{2}}\right) \cos\left(\frac{9 \pi y}{50}\right)\right)
	\end{split}}
\end{equation}
 is an approximation of a non-uniform stationary solution of \eqref{FunctionResp}-\eqref{NewmanCond}. Using \eqref{appr3} as initial data we see the evolution in time and the convergence of the Midpoint Method to an exact non-uniform stationary solution of system \eqref{FunctionResp}-\eqref{NewmanCond}, see \cite[\href{https://tinyurl.com/8ywkptst}{\blue Interactive Simulation 3}]{Walker et all}.
  \begin{figure}[h!]
 	\begin{center}
 	\includegraphics[width=6cm]{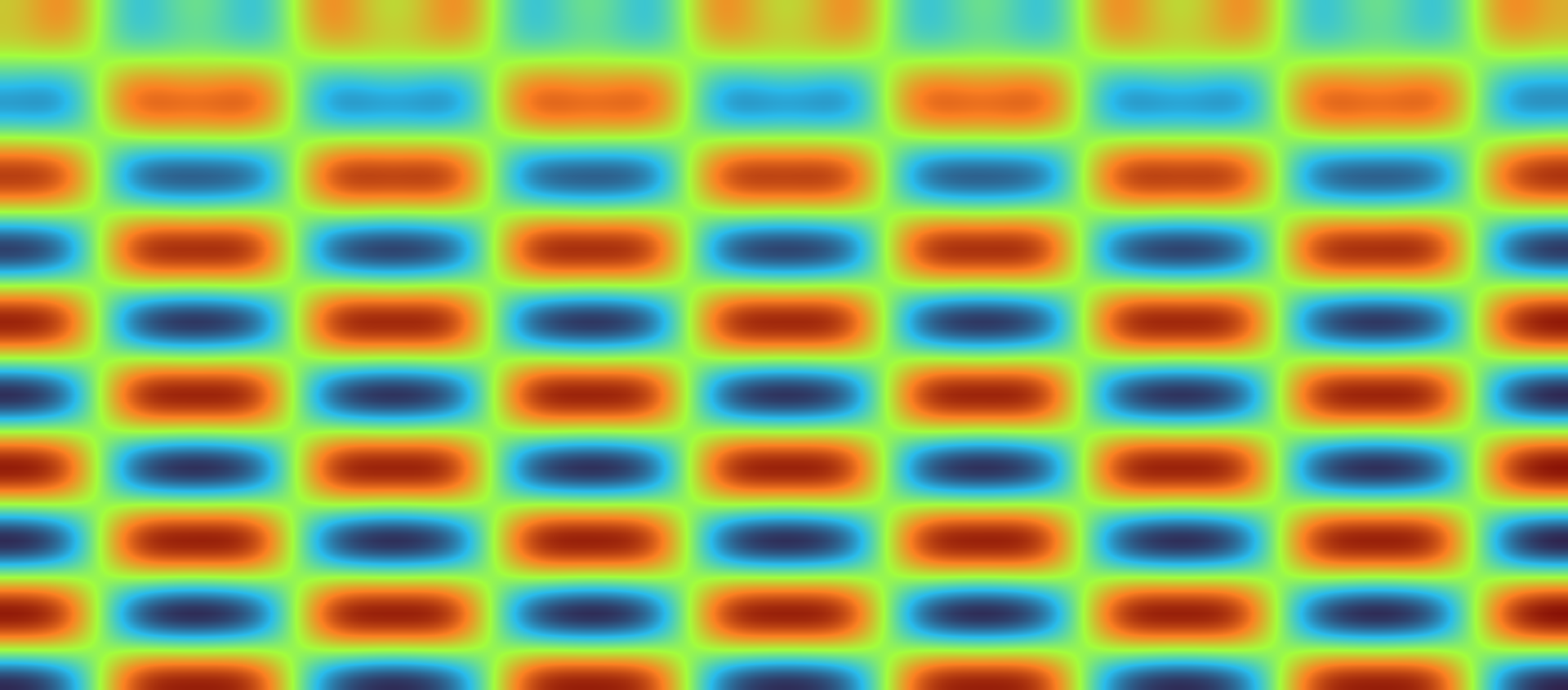}	
 	\includegraphics[width=7cm]{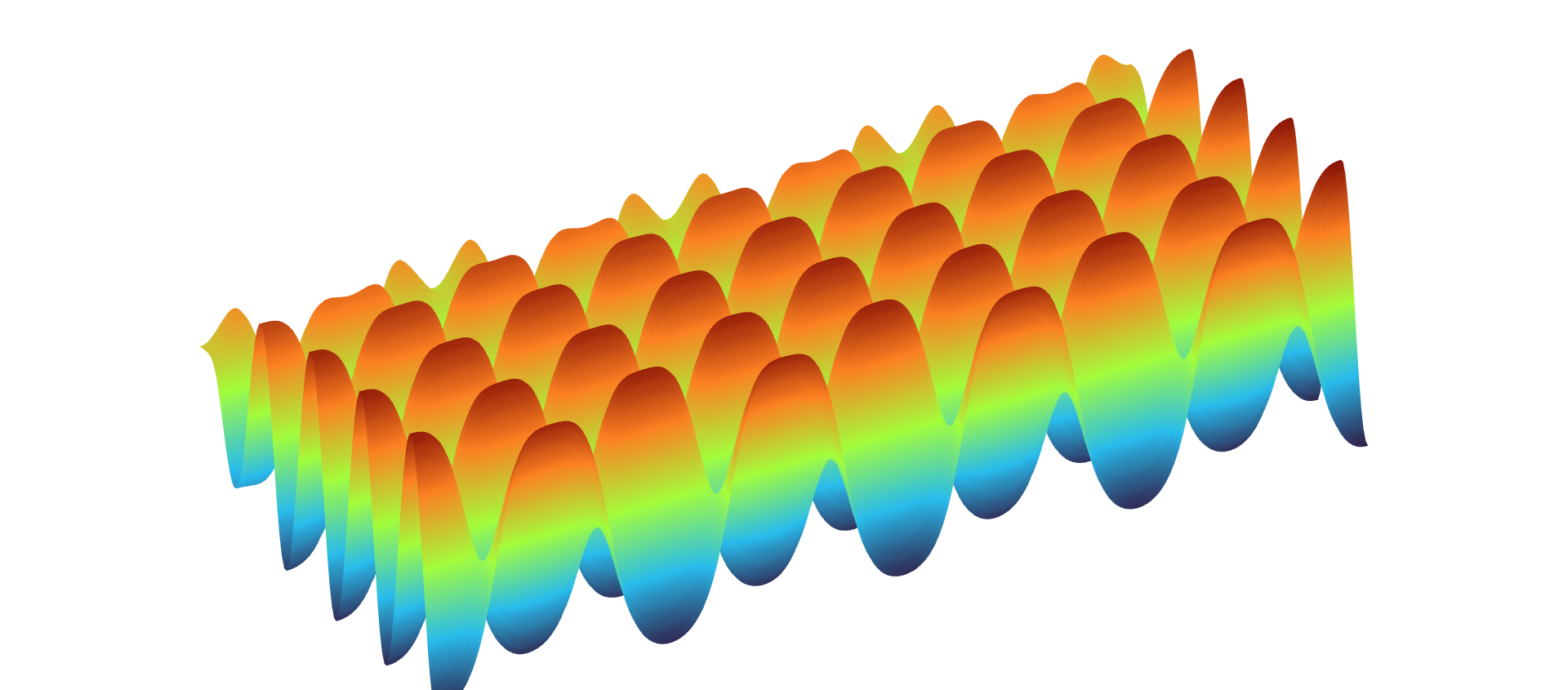}
 		\caption{The population density of prey $0.493 \leq u(x) \leq 0.507.$}\label{figura12}
 	\end{center}
 \end{figure}
  \begin{figure}[h!]
 	\begin{center}
 		\includegraphics[width=6cm]{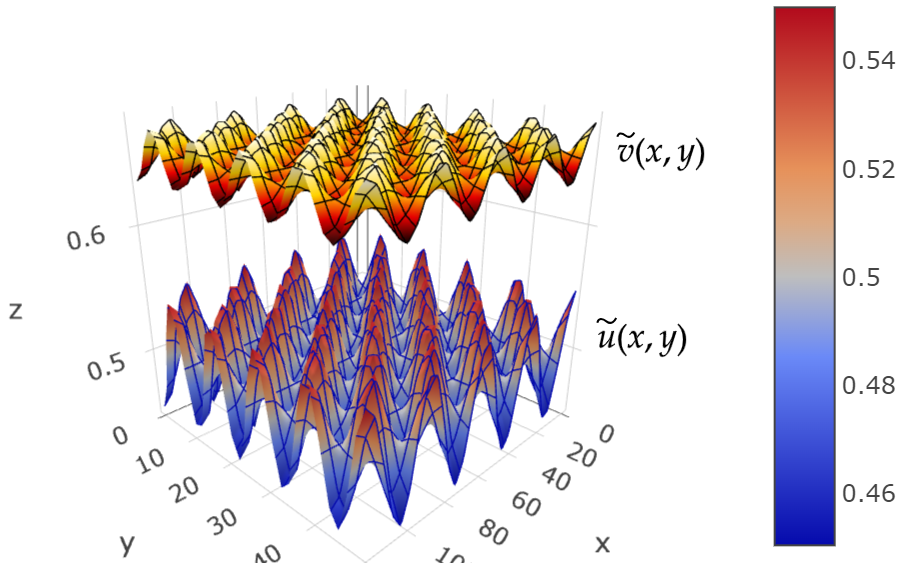}
 		\caption{Projection of the population density approximations of predator $\tilde{v}$ and the prey $\tilde{u}$ on the $xyz-$plane.}\label{figura13}
 	\end{center}
 \end{figure}
\begin{rem}
	As stated in \cite{Bunow et al}, changes in domain shape, produce changes in the bifurcated patterns. For instance, if we choose a square ($L_{1}=L_{2}=100$) in \eqref{unidlapl4} and consider the simple eigenvalue $\tilde{\lambda}_{9,9}\approx 0.1599,$ and $d_{1}=1.0733,$ Then at $d_{2}^{\ast}=232.4321$ the uniform steady-state solution  $w_{0}=(0.4998, 0.6532) \in \mathcal{C}_{1} \cap \mathcal{C}_{2} \cap \mathbb{R}^{2}_{+}$ of system \eqref{FunctionResp}-\eqref{NewmanCond} undergoes a Turing bifurcation, see \cite[\href{https://tinyurl.com/4xz2pn4v}{\blue Interactive Simulation 4}]{Walker et all} and Figure \ref{figura15}.
\end{rem}
 \section{Concluding Remarks}\label{CR} 
 In this work, we have  shown two different criteria to prove the existence of patterns for reaction-diffusion models of two components in the space of continuous functions (see \eqref{Gensyst3}). Both criteria reveal the importance of knowing the simple eigenvalues of the Laplacian operator on bounded, open and connected domains. We have applied these results to obtain sufficient conditions on the parameters present in a particular nonlinear predator-prey model (see \eqref{FunctionResp}-\eqref{NewmanCond}) to prove that it can exhibit stable spatially heterogeneous solutions which arise from Turing bifurcations. We observe that for such a particular reaction-diffusion model, a Turing Bifurcation can not occur for a large diffusive coefficient of prey, nevertheless, the diffusive coefficient of predator can be large enough. Furthermore, the gap between the maximal mortality $\delta$ and the mortality at low density $\gamma$ can be large as well as the conversion rate $b.$ This analysis shows how a reaction-diffusion predator-prey model with variable mortality and Hollyn's type II functional response can stably regulate its growth around either spatially homogeneous or heterogeneous solutions through a Turing instability mechanism.
 
 From this article, a question arises: It would be possible to show the pattern formation for reaction-diffusion systems of two components  by using a not simple eigenvalue (eigenvalue with multiplicity greater than or equal to two) in the hypotheses of Criteria I or II? The answer to  this question and   similar criterion for  reaction-diffusion models of three components will be discussed in a forthcoming paper. This work is in progress.

 \begin{appendix}\label{ppenA}
 \section{}
 
 In this subsection, we introduce the concept of excitable matrix and recall conditions under which  real $k\times k,$  matrices are excitable. We refer the reader to \cite{Cross, Hershkowitz, Satnoianu et al 1, Satnoianu et al 2, Hadeler and Ruan} for more details. We denote by $\mathcal{M}_{k\times k}(\mathbb{R})$ the space of real $k\times k$ matrices with real coefficients.
 
 \begin{defn}[Matrix Stability. See \cite{Hadeler and Ruan, Cross}]\label{stabmatrx}
 	Let $A \in \mathcal{M}_{k\times k}(\mathbb{R}).$ 
 	\begin{itemize}
 		\item[$i)$] $A$ is said to be stable if all eigenvalues of $A$ are located in the open left half-plane of the complex plane.
 		
 		\item[$ii)$] $A$ is said to be strongly stable (with respect to diffusion) if $A-M$ is stable for any nonnegative definite diagonal real $k\times k$ matrix $M.$
 		
 		\item[$iii)$] $A$ is said to be excitable (with respect to diffusion) if $A$ is stable but not strongly stable.		
 	\end{itemize}
 \end{defn}

 Of course, a strongly stable matrix is also stable. Also, for an excitable matrix $A$ there is always a choice of a non-negative definite diagonal real matrix $M$ such that $A-M$ is unstable. 
 
 \begin{defn}\label{minor}
 	Let $A \in \mathcal{M}_{k\times k}(\mathbb{R}).$ 
 	For any subset $1\leq i_{1}<i_{2} \cdots <i_{j}\leq k$ of the integers $1,2,...,k,$ the square principal submatrix $A_{i_{1},i_{2},...,i_{j}}$ of $A$ is obtained by taking exactly the rows and the columns of indices $i_{1},i_{2},...,i_{j}.$ The corresponding principal minor (of order $j$ for $A$) is $M_{i_{1},i_{2},...,i_{j}}:=\text{det}(A_{i_{1},i_{2},...,i_{j}}).$ The signed principal minors of $A$ are the quantities $(-1)^{j}M_{i_{1},i_{2},...,i_{j}}.$ The minors $M_{i}$ are written simply $a_{ii}.$\\
 	
 	We said that $A$ is s-stable if, for any minor 
 	$M_{i_{1},i_{2},...,i_{j}}$ of order $j$ ($1\leq j \leq k $), we have $$\text{sgn}(M_{i_{1},i_{2},...,i_{j}})=(-1)^{j}.$$

 \end{defn}
 We denote by $\mathcal{P}$ the class of matrices whose signed principal minors are all positive and by $\mathcal{P}_{0}^{+}$ the class of matrices whose signed principal minors are all non-negative, with at least one of each order positive. A complete characterization of strongly stable matrices, in terms of inequalities on their minors, has been given by Cross \cite{Cross} for $k=2$ and $k=3$ using the Routh-Hurwitz condition.
 
 \begin{thm}[See Theorem $3$ in \cite{Cross}]\label{strsta2x2}
 	Let $A=(a_{ij}),$ $1\leq i, j\leq 2$ be a matrix in $\mathcal{M}_{2\times 2}(\mathbb{R}).$ Then, $A$ is strongly stable if and only if $A \in \mathcal{P}_{0}^{+}.$ Therefore, $A$ is strongly stable if the following conditions holds:
 	\begin{equation}\label{strsta2x21}
 		\text{Tr}(A)=a_{11}+a_{22}<0,
 	\end{equation}	
 	\begin{equation}\label{strsta2x22}
 		\text{det}(A)=a_{11}\cdot a_{22}-a_{21}\cdot a_{12}>0,
 	\end{equation}
 	\begin{equation}\label{strsta2x23}
 		a_{11}\leq 0,\;\;\;a_{22}\leq 0.
 	\end{equation} 	 
 	
 	The matrix $A$ is excitable if occur \eqref{strsta2x21}, \eqref{strsta2x22} and 
 	\begin{equation}\label{excita2x2}
 		a_{11}> 0\;\;\text{or}\;\;a_{22}> 0.
 	\end{equation} 	
 \end{thm}
 
 \begin{rem}\label{excita2x21}
 	Hence, there are only four sign arrangements for matrices $A=(a_{ij}),$ $1\leq i, j\leq 2$  in $\mathcal{M}_{2\times 2}(\mathbb{R})$ to be excitable:
 	$$ \left(
 	\begin{array}{cc}
 		+&- \\
 		+&-\\
 	\end{array}
 	\right);
 	\left(
 	\begin{array}{cc}
 		+&+ \\
 		-&-\\
 	\end{array}
 	\right);
 	\left(
 	\begin{array}{cc}
 		-&+ \\
 		-&+\\
 	\end{array}
 	\right);
 	\left(
 	\begin{array}{cc}
 		-&- \\
 		+&+\\
 	\end{array}
 	\right). $$	
 \end{rem}
 
 The necessary and sufficient conditions for matrices of order $3$ are stated in the following result:

 \begin{thm}[See Theorem $4$ in \cite{Cross}]\label{strsta3x3}
 	Let $A=(a_{ij}),$ $1\leq i, j\leq 3$ be a matrix in $\mathcal{M}_{3\times 3}(\mathbb{R}).$ Then, $A$ is strongly stable if and only if $A \in \mathcal{P}_{0}^{+}$ and $A$ is stable. Hence, $A$ is strongly stable if and only if the following conditions holds:
 	\begin{equation}\label{strsta3x31}
 		-c_{1}:=\text{Tr}(A)
 		=a_{11}+a_{22}+a_{33}<0,
 	\end{equation}	
 	\begin{equation}\label{strsta3x32}
 		-c_{3}:=\text{det}(A)<0,
 	\end{equation}
 	\begin{equation}\label{strsta3x34}
 		c_{1}\cdot c_{2}-c_{3}>0,
 	\end{equation} 
 	\begin{equation}\label{strsta3x33}
 		a_{11}\leq 0,\;\;\;a_{22}\leq 0, \;\;\text{and}\;\;a_{33}\leq 0,
 	\end{equation} 	 
 	\begin{equation}\label{strsta3x35}
 		M_{12}=a_{11} a_{22}-a_{12} a_{21}\geq 0,\;\;\;M_{13}=a_{11} a_{33}-a_{13} a_{31}\geq 0, \;\;\text{and}\;\;M_{23}=a_{22} a_{33}-a_{23} a_{32}\geq 0,
 	\end{equation}
 	where $c_{2}:=M_{12}+M_{13}+M_{23}.$\\
 	
 	The matrix $A$ is excitable if occur \eqref{strsta3x31}, \eqref{strsta3x32}, \eqref{strsta3x34} and 
 	\begin{equation}\label{strsta3x36}
 		a_{11}> 0,\;\;\text{or}\;\;a_{22}> 0, \;\;\text{or}\;\;a_{33}> 0,
 	\end{equation} 	 
 	or
 	\begin{equation}\label{strsta3x37}
 		M_{12}=a_{11} a_{22}-a_{12} a_{21}< 0,\;\text{or}\;M_{13}=a_{11} a_{33}-a_{13} a_{31}< 0, \;\text{or}\; M_{23}=a_{22} a_{33}-a_{23} a_{32}< 0.
 	\end{equation}
 \end{thm}
 
 \begin{rem}\label{excita3x31}
 	For matrices of order $3,$ even knowing the signs of the matrix $A,$ we must request extra conditions on its coefficients so that $A$ is excitable. 
 \end{rem}
 
 Characterizing all strongly stables matrices for $k\geq 4$ is an open problem (see \cite{Satnoianu et al 2}). Nevertheless, we have the following result establishing  necessary conditions for a matrix to be excitable.
 
 \begin{thm}[See Theorem 2 in \cite{Satnoianu et al 1}]\label{neccond}
 	Let $A \in \mathcal{M}_{k\times k}(\mathbb{R}).$ If $A$ is strongly stable then $A$ is s-stable. If the matrix $A$ is stable but not s-stable, then $A$ is excitable .
 \end{thm}
 
 \begin{rem}\label{necexcitable}
 	Theorem \eqref{neccond} says that a sufficient condition for a stable matrix	$A \in \mathcal{M}_{k\times k}(\mathbb{R})$ to be excitable, is that there exists a principal minor $M_{i_{1},i_{2},...,i_{j}}$ of order $j$ for $A$ ($1\leq j < k $) such that
 	$\text{sgn}(M_{i_{1},i_{2},...,i_{j}})\neq (-1)^{j}.$ 
 \end{rem}

 The following nomenclature serves to distinguish the different possible cases.
 \begin{defn}\label{typeexcmat}
 	Let $A \in \mathcal{M}_{k\times k}(\mathbb{R})$ be an excitable matrix.
 	\begin{itemize}
 		\item[$i)$] We said that $A$ is an excitable matrix of the first type if $M_{i}=a_{ii}>0$ for some $1\leq i\leq k.$
 		\item[$ii)$] We said that $A$ is an excitable matrix of the second type if there are $1\leq i,j\leq k$ such that $M_{ij}<0$ and $A$ is not of the first type.
 		\item[$iii)$]	Inductively, we said that $A$ is an excitable matrix of type $j,$ with $j<k,$ if there are indices $1\leq i_{1},i_{2},...,i_{j}< k$ such that $\text{sgn}(M_{i_{1}j_{2}...i_{j}})=(-1)^{j+1},$ and $A$ is not being associated with an excitable matrix of any previous type.	 
 	\end{itemize}
 \end{defn}
 
 \end{appendix}


\subsection*{Acknowledgment}
Vielma-Leal F. acknowledges that this was a project supported by the Competition for Research Regular Projects, year 2022, code LPR22-02, Universidad Tecnológica Metropolitana.  Del Rio Palma M. was supported by the Conselho Nacional de Desenvolvimento Científico e Tecnológico (CNPq) and the FAPERJ-Fundação Carlos Chagas Filho de Amparo à Pesquisa do Estado do Rio de Janeiro Processo SEI 260003/014835/2023 and CNPq 151052/2023-9 and  
Montenegro-Concha M. was partially supported by STEMDEPID, year 2022, code  AMSUD 220002, MathAmsud 2022.


\end{document}